\documentclass{amsart}

\input epsf

\title[The Shape of a Typical Boxed Plane Partition]{The Shape of a
Typical Boxed Plane Partition}
\author{Henry Cohn}
\address{Department of Mathematics\\
Harvard University\\
Cambridge, MA 02138}
\email{cohn@math.harvard.edu}
\author{Michael Larsen}
\address{Department of Mathematics\\
Indiana University\\
Bloomington, IN 47405}
\email{larsen@math.missouri.edu}
\author{James Propp}
\address{Department of Mathematics\\
University of Wisconsin\\
Madison, WI 53706}
\email{propp@math.wisc.edu}

\date{May 28, 2002;  this version of the paper is updated from
the published version (New York Journal of Mathematics {\bf 4}
(1998), 137--165) by correcting four inconsequential typos
and adding more bibliographic information for several references.}

\thanks{Cohn was supported by an NSF Graduate Research Fellowship.
Larsen was supported by NSF grant DMS94-00833 and by the Sloan
Foundation.  Propp was supported by NSA grant MDA904-92-H-3060 and NSF
grant DMS92-06374, and by a grant from the MIT Class of 1922.}



\theoremstyle{plain}
\newtheorem{thm}{Theorem}[section]

\newtheorem{prop}[thm]{Proposition}

\theoremstyle{definition}
\newtheorem{conj}[thm]{Conjecture}
\newtheorem{quest}[thm]{Open Question}

\numberwithin{equation}{section}

\newcommand{\thmref}[1]{Theorem~\ref{#1}}
\newcommand{\propref}[1]{Proposition~\ref{#1}}

\newcommand{\conjref}[1]{Conjecture~\ref{#1}}
\newcommand{\questref}[1]{Open Question~\ref{#1}}

\newcommand{\R} {{\mathbb R}}

\newcommand{\V} {{\mathcal V}}
\newcommand{\A} {{\mathbb A}}
\newcommand{\Admis} {{\mathcal F}}
\newcommand{\Hf} {{\mathcal H}}
\newcommand{\Pl} {{\mathcal P}}
\newcommand{\ProjP} {{\mathbb P}}

\newcommand{\sumL} {\!\!\mathop{{\,\,\,}{\sum}^{_L}}}

\newcommand{\Real} {\textup{Re}\,}
\newcommand{\Imag} {\textup{Im}\,}

\newcommand{\Philim} {F}

\newcommand{\sgn} {\textup{sgn}\,}
\newcommand{\Res} {\textup{Res}}
\newcommand{\intr} {\int_{-\infty}^\infty}
\newcommand{\inprod}[2] {{\langle #1, #2 \rangle}}

\newcommand{\nleft}{n_L}
\newcommand{\nright}{n_R}
\newcommand{\nmiddle}{n}
\newcommand{\nbig}{N}
\newcommand{\rhL}{\rho_L}
\newcommand{\rhR}{\rho_R}
\newcommand{\rhoLR}{\rho}
\newcommand{\scalefactor}{\sigma}
\newcommand{\oa}{a'}
\newcommand{\ob}{b'}
\newcommand{\oc}{c'}

\newcommand{\SL}{\mbox{SL}}

\newcommand{\abstracttext}{ Using a calculus of variations approach,
we determine the shape of a typical plane partition in a large box
(i.e., a plane partition chosen at random according to the uniform
distribution on all plane partitions whose solid Young diagrams fit
inside the box).  Equivalently, we describe the distribution of the
three different orientations of lozenges in a random lozenge tiling of
a large hexagon.  We prove a generalization of the classical formula
of MacMahon for the number of plane partitions in a box; for each of
the possible ways in which the tilings of a region can behave when
restricted to certain lines, our formula tells the number of tilings
that behave in that way.  When we take a suitable limit, this formula
gives us a functional which we must maximize to determine the
asymptotic behavior of a plane partition in a box.  Once the
variational problem has been set up, we analyze it using a
modification of the methods employed by Logan and Shepp and by Vershik
and Kerov in their studies of random Young tableaux.}

\begin{document}

\begin{abstract}
\abstracttext
\end{abstract}

\maketitle

\tableofcontents

\section{Introduction}
\label{sec-intro}

In this paper we will show that almost all plane partitions that are
constrained to lie within an $a \times b \times c$ box have a certain
approximate shape, if $a$, $b$, and $c$ are large; moreover, this limiting
shape depends only on the mutual ratios of $a$, $b$, and $c$.  Our proof will
make use of the equivalence between such plane partitions and tilings of
hexagons by rhombuses.

Recall that plane partitions are a two-dimensional generalization of ordinary
partitions.  A {\it plane partition} $\pi$ is a collection of non-negative
integers $\pi_{x,y}$ indexed by ordered pairs of non-negative integers $x$ and
$y$, such that only finitely many of the integers $\pi_{x,y}$ are non-zero,
and for all $x$ and $y$ we have $\pi_{x+1,y} \le \pi_{x,y}$ and $\pi_{x,y+1}
\le \pi_{x,y}$.  A more symmetrical way of looking at a plane partition is to
examine the union of the unit cubes $[i,i+1] \times [j,j+1] \times [k,k+1]$
with $i$, $j$, and $k$ non-negative integers satisfying $0 \le k < \pi_{i,j}$.
This region is called the {\it solid Young diagram} associated with the plane
partition, and its volume is the sum of the $\pi_{x,y}$'s.

\begin{figure}
\begin{center}
\leavevmode
\epsfbox[0 0 120 104]{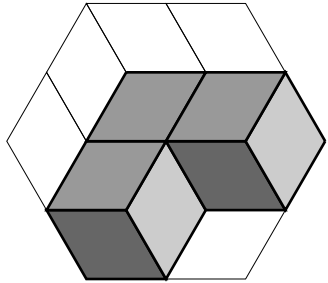}
\end{center}
\caption{The tiling corresponding to a plane partition.}
\label{fig-shadedhex}
\end{figure}

We say that a plane partition $\pi$ {\it fits within an $a \times b \times c$
box} if its solid Young diagram fits inside the rectangular box $[0,a] \times
[0,b] \times [0,c]$, or equivalently, if $\pi_{x,y} \le c$ for all $x$ and
$y$, and $\pi_{x,y}=0$ whenever $x \ge a$ or $y \ge b$; we call such a plane
partition a {\it boxed plane partition}.  Plane partitions in an $a \times b
\times c$ box are in one-to-one correspondence with tilings of an equi-angular
hexagon of side lengths $a,b,c,a,b,c$ by rhombuses whose sides have length $1$
and whose angles measure $\pi/3$ and $2\pi/3$.  It is not hard to write down a
bijection between the plane partitions and the tilings, but the correspondence
is best understood informally, as follows.  The faces of the unit cubes that
constitute the solid Young diagram are unit squares.
Imagine now that we augment the solid Young diagram by adjoining the three
``lower walls'' of the $a \times b \times c$ box that contains it (namely
$\{0\} \times [0,b] \times [0,c]$, $[0,a] \times \{0\} \times [0,c]$, and
$[0,a] \times [0,b] \times \{0\}$); imagine as well that each of these walls
is divided into unit squares.  If we look at this augmented Young diagram
{}from a point on the line $x=y=z$, certain of the unit squares are visible
(that is, unobstructed by cubes).  These unit squares form a surface whose
boundary is the non-planar hexagon whose vertices are the points $(a,0,c)$,
$(0,0,c)$, $(0,b,c)$, $(0,b,0)$, $(a,b,0)$, $(a,0,0)$, and $(a,0,c)$,
respectively.  Moreover, these same unit squares, projected onto the plane
$x+y+z=0$ and scaled, become rhombuses which tile the aforementioned
planar hexagon.
For example, the plane partition $\pi$
in a $2 \times 2 \times 2$ box defined by $\pi_{0,0} = \pi_{0,1} = \pi_{1,0} =
1$ and $\pi_{1,1}=0$ corresponds to the tiling in Figure~\ref{fig-shadedhex},
where the points $(2,0,0)$, $(0,2,0)$, and $(0,0,2)$ are at the lower left,
extreme right, and upper left corners of the hexagon.  (The shading is meant
as an aid for three-dimensional visualization, and is not necessary
mathematically.  The unshaded rhombuses represent part of the walls.)

We will use the term {\it $a,b,c$ hexagon} to refer to an equi-angular
hexagon of side lengths $a,b,c,a,b,c$ (in clockwise order, with the
horizontal sides having length $b$), and the term {\it lozenge} to
refer to a unit rhombus with angles of $\pi/3$ and $2\pi/3$.  We will
focus, without loss of generality, on those lozenges whose major axes
are vertical, which we call {\it vertical lozenges}.  Although our
method in this article is to reduce facts about plane partitions to
facts about tilings, one can also go in the reverse direction.  For
example, one can see from the three-dimensional picture that in every
lozenge tiling of an $a,b,c$ hexagon, the number of vertical lozenges
is exactly $ac$ (with similar formulas for the other two orientations
of lozenge); see \cite{DT} for further discussion.

A classical formula of MacMahon \cite{M} asserts that there are
exactly
$$
\prod_{i=0}^{a-1} \prod_{j=0}^{b-1} \prod_{k=0}^{c-1}
\frac{i+j+k+2}{i+j+k+1}
$$
plane partitions that fit in an $a \times b \times c$ box, or
(equivalently) lozenge tilings of an $a,b,c$ hexagon.  In
\thmref{countonline} of this article, we give a generalization that
counts lozenge tilings with prescribed behavior on a given horizontal
line.

Using \thmref{countonline}, we will determine the shape of a typical
plane partition in an $a \times b \times c$ box (\thmref{heightfun}).
Specifically, that theorem implies that if $a,b,c$ are large, then the
solid Young diagram of a random plane
partition in an $a \times b \times c$ box is almost certain to differ
{}from a particular, ``typical'' shape by an amount that is negligible
compared to $abc$ (the total volume of the box).  Equivalently, the
visible squares in the augmented Young diagram of the random boxed
plane partition
form a surface whose maximum deviation from a
particular, typical surface is almost certain to be
$o(\min(a,b,c))$.  Moreover, scaling $a,b,c$ by some factor causes the
typical shape of this bounding surface to scale by the same factor.

Before we say what the true state of affairs is, we invite the reader
to come up with a guess for what this typical shape should be.  One
natural way to arrive at a guess is to consider the analogous problem
for ordinary (rather than plane) partitions.  If one considers all
ordinary partitions that fit inside an $a \times b$ rectangle (in the
sense that their Young diagrams fit inside $[0,a] \times [0,b]$), then
it is not hard to show that the typical boundary of the diagram is the
line $x/a+y/b=1$; that is, almost all such partitions have roughly
triangular Young diagrams.  (One way to prove this is to apply
Stirling's formula to binomial coefficients and employ direct
counting; another is to use probabilistic methods, aided by a
verification
that if we look at the boundary of the Young diagram of the partition
as a lattice path, then the directions of different steps in the path
are negatively correlated.)
It therefore might seem plausible that the typical bounding
surface for plane partitions would be a plane (except where it coincided with
the sides of the box), as when a box is tilted on its corner and half-filled
with fluid.  However, \thmref{mainthm} shows that that is in fact not the
case.

\begin{figure}
\begin{center}
\leavevmode
\epsfbox[0 0 288 250]{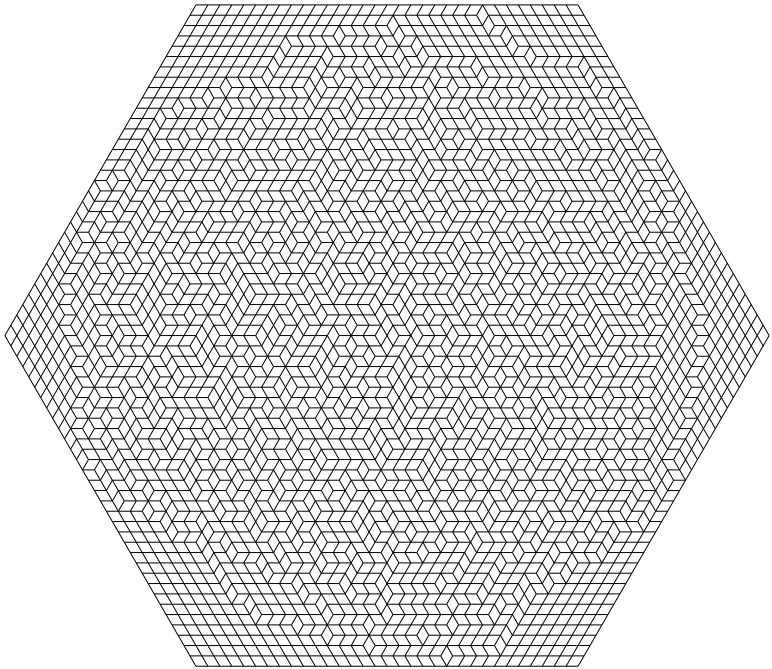}
\end{center}
\caption{A random lozenge tiling of a $32,32,32$ hexagon.}
\label{fig-random}
\end{figure}

To see what a typical boxed plane partition {\it does} look like, see
Figure~\ref{fig-random}.  This tiling was generated using the methods
of \cite{PW} and is truly random, to the extent that pseudo-random
number generators can be trusted.  Notice that near the corners of the
hexagon, the lozenges are aligned with each other, while in the
middle, lozenges of different orientations are mixed together.  If the
bounding surface of the Young diagram tended to be flat, then the
central zone of mixed orientations would be an inscribed hexagon, and
the densities of the three orientations of tiles would change
discontinuously as one crossed into this central zone.  In fact, what
one observes is that the central zone is roughly circular, and that
the tile densities vary continuously except near the midpoints of the
sides of the original hexagon.

One can in theory use our results to obtain an explicit formula for
the typical shape of the bounding surface, in which one specifies the
distance from a point $P$ on the surface to its image $P'$ under
projection onto the $x+y+z=0$ plane, as a function of $P'$; however,
the integral that turns up is quite messy (albeit evaluable in closed
form), with the result that the explicit formula is extremely lengthy
and unenlightening.  Nevertheless, we can and do give a fairly simple
formula for the tilt of the tangent plane at $P$ as a function of the
projection $P'$, which would allow one to reconstruct the surface
itself.  In view of the correspondence between plane partitions and
tilings, specifying the tilt of the tangent plane is equivalent to
specifying the local densities of the three different orientations of
lozenges for random tilings of an $a,b,c$ hexagon.

Our result on local densities has as a particular corollary the
assertion that, in an asymptotic sense, the zone of mixed orientation
(defined as the region in which all three orientations of lozenge
occur with positive density) is precisely the interior of the ellipse
inscribed in the hexagon.  This behavior is analogous to what has been
proved concerning domino tilings of regions called Aztec diamonds (see
\cite{JPS} and \cite{CEP}); these are roughly square regions, and if
one tiles them randomly with dominos ($1 \times 2$ rectangles), then
the zone of mixed orientation tends in the limit to the inscribed
disk.  However, the known results for Aztec diamonds are stronger than
the corresponding best known results for hexagons (see \conjref{aec}
in Section~\ref{sec-conj}).

To state our main theorem, we begin by setting up normalized
coordinates.  Suppose that we are dealing with an $a,b,c$ hexagon, so
that the side lengths are $a,b,c,a,b,c$ in clockwise order (with the
sides of length $b$ horizontal).  We let $a,b,c$ tend to infinity
together, in such a way that the three-term ratio $a:b:c$ (i.e.,
element of the
projective plane $\ProjP^2(\R)$) converges to $\alpha:\beta:\gamma$
for some fixed positive numbers $\alpha$, $\beta$, $\gamma$.  Say
$a+b+c=\scalefactor(\alpha+\beta+\gamma)$ for some scaling factor
$\scalefactor$.  Then we choose re-scaled coordinates for the $a,b,c$
hexagon so that its sides are $\oa=a/\scalefactor$,
$\ob=b/\scalefactor$, and $\oc=c/\scalefactor$ (which by the
hypothesis of our main theorem will be required to converge to
$\alpha$, $\beta$, and $\gamma$ as $a$, $b$, and $c$ get large).  Note
that $\oa,\ob,\oc$ in general are not integers.  The origin of our
coordinate system lies at the center of the hexagon.  One can check
that the sides of the hexagon lie on the lines
$y=\frac{\sqrt{3}}{2}(2x+\ob+\oc)$, $y=\frac{\sqrt{3}}{4}(\oa+\oc)$,
$y=\frac{\sqrt{3}}{2}(-2x+\oa+\ob)$,
$y=\frac{\sqrt{3}}{2}(2x-\ob-\oc)$, $y=\frac{\sqrt{3}}{4}(-\oa-\oc)$,
and $y=\frac{\sqrt{3}}{2}(-2x-\oa-\ob)$, and that the inscribed
ellipse is described by
$$
3(\oa+\oc)^2x^2 - 2\sqrt{3}(\oa+2\ob+\oc)(\oa-\oc)xy+
((\oa+2\ob+\oc)^2-4\oa\oc)y^2
= 3\oa\ob\oc(\oa+\ob+\oc).
$$
Define $E_{\alpha,\beta,\gamma}(x,y)$ to be the polynomial
$$
3\alpha\beta\gamma(\alpha+\beta+\gamma) - (3(\alpha+\gamma)^2x^2
- 2\sqrt{3}(\alpha+2\beta+\gamma)(\alpha-\gamma)xy
+ ((\alpha+2\beta+\gamma)^2-4\alpha\gamma)y^2),
$$
whose zero-set is the ellipse inscribed in the $\alpha,\beta,\gamma$
hexagon.
Also define
$Q_{\alpha,\beta,\gamma}(x,y)$ to be the polynomial
$$
\frac{\sqrt{3}}{2}\left(
\frac{4}{3}y^2-4x^2+\beta^2+\alpha\beta+\beta\gamma-\alpha\gamma\right),
$$
which will be useful shortly.

There are six points at which the ellipse inscribed in the
$\alpha,\beta,\gamma$ hexagon meets the boundary of the hexagon.  The
four that occur on sides of length $\alpha$ or $\gamma$ will be called
{\it singular points}, for reasons that will be clear shortly.
(Recall that we have already broken the underlying symmetry of the
situation by agreeing to focus on vertical lozenges.)
The points of the hexagon that lie outside the inscribed ellipse form six
connected components.  Let $R_1$ be the closure of the union of the
two components containing the leftmost and rightmost points of the
hexagon, minus the four singular points, and let $R_0$ be the closure
of the union of the other four components, again minus the four
singular points.

Finally, define the (normalized) coordinates of a vertical lozenge to
be the (normalized) coordinates of its center.  Then the following
theorem holds:

\begin{thm}
\label{mainthm}
Let $U$ be the interior of a smooth simple closed curve contained inside
the $\alpha,\beta,\gamma$ hexagon,
with $\alpha,\beta,\gamma > 0$.
In a random tiling of an $a,b,c$ hexagon, as $a,b,c \rightarrow \infty$
with $a:b:c \rightarrow \alpha:\beta:\gamma$, the expected
number of vertical lozenges whose normalized coordinates lie in $U$
is $ab+bc+ac$ times
$$
\frac{1}{A}
\iint_{\!U} \Pl_{\alpha,\beta,\gamma}(x,y) \, dx \, dy \, + \, o(1),
$$
where $A = (\alpha\beta+\beta\gamma+\alpha\gamma)\sqrt{3}/2$
is the area of the
$\alpha,\beta,\gamma$ hexagon and
$$
\Pl_{\alpha,\beta,\gamma}(x,y) = \begin{cases}
\displaystyle
\frac{1}{\pi}
\cot^{-1}\!\left(\frac{Q_{\alpha,\beta,\gamma}(x,y)}
{\sqrt{E_{\alpha,\beta,\gamma}(x,y)}}\right) &
\textup{if $E_{\alpha,\beta,\gamma}(x,y) > 0$,} \\
0 & \textup{if $(x,y) \in R_0$, and}\\
1 & \textup{if $(x,y) \in R_1$}.\\
\end{cases}
$$
\end{thm}

In fact, our proof will give an even stronger version of
\thmref{mainthm}, in which $U$ is a horizontal line segment rather
than an open set (and the double integral is replaced by a single
integral).  Since we can derive the expected number of vertical
lozenges whose normalized coordinates lie in $U$ by integrating over
all horizontal cross sections (provided that the error term is uniform
for all cross sections, as we will prove in Section~\ref{sec-analysis}),
this variant of the claim implies the
one stated above, though it is not obviously implied by it.  We have
stated the result in terms of open sets because that formulation seems
more natural; the proof just happens to tell us more.

The intuition behind \thmref{mainthm} is that
$\Pl_{\alpha,\beta,\gamma}(x,y)$ gives the density of vertical
lozenges in the normalized vicinity of $(x,y)$; the factor of
$ab+bc+ac$ arises simply because there are that many lozenges in a
tiling of an $a,b,c$ hexagon.  One might be tempted to go farther and
think of $\Pl_{\alpha,\beta,\gamma}(x,y)$ as the asymptotic
probability that a random tiling of the $a,b,c$ hexagon will have a
vertical lozenge at any particular location in the normalized vicinity
of $(x,y)$ (see \conjref{cep} in Section~\ref{sec-conj}); however, we
cannot justify this interpretation rigorously, because it is
conceivable that there are small-scale fluctuations in the
probabilities that disappear in the double integral.  (In
Subsection~1.3 of \cite{CEP}, it is shown that the analogous
probabilities for random domino tilings do exhibit such fluctuations,
although the fluctuations disappear if one distinguishes between four
classes of dominos, rather than just horizontal and vertical dominos.)

The formula for $\Pl_{\alpha,\beta,\gamma}$ is more natural than it might
appear.  Examination of random tilings such as the one shown in
Figure~\ref{fig-random} leads one to conjecture that the region in which all
three orientations of lozenges occur with positive density is (asymptotically)
the interior of the ellipse inscribed in the hexagon, and the known fact that
an analogous claim holds for random domino tilings of Aztec diamonds (see
\cite{CEP}) lends further support to this hypothesis.  This
leads one to predict that $\Pl_{\alpha,\beta,\gamma}$ will be $0$ in $R_0$ and
$1$ in $R_1$, and strictly between $0$ and $1$ in the interior of the ellipse.
Comparison with the analogous theorem for domino tilings (Theorem~1 of
\cite{CEP}) suggests that within the inscribed ellipse,
$\Pl_{\alpha,\beta,\gamma}(x,y)$ should be of the form
$$
\frac{1}{\pi}
\cot^{-1}\left(\frac{Q(x,y)}{\sqrt{E_{\alpha,\beta,\gamma}(x,y)}}\right)
$$
for some quadratic polynomial $Q(x,y)$, and the only problem is
actually finding $Q(x,y)$ in terms of $\alpha,\beta,\gamma$.  A simple
description of the polynomial $Q(x,y)$ that actually works is that its
zero-set is the unique hyperbola whose asymptotes are parallel to the
non-horizontal sides of the hexagon and which goes through the four
points on the boundary of $R_1$ where the inscribed ellipse is tangent
to the hexagon.  This determines $Q(x,y)$ up to a constant factor, and
that constant factor is determined (at least in theory---in practice
the calculations would be cumbersome) by requiring that the average of
$\Pl_{\alpha,\beta,\gamma}(x,y)$ over the entire hexagon must be
$\alpha\gamma/(\alpha\beta+\beta\gamma+\alpha\gamma)$.

\begin{figure}
\begin{center}
\leavevmode
\epsfbox[-8 -8 130 114]{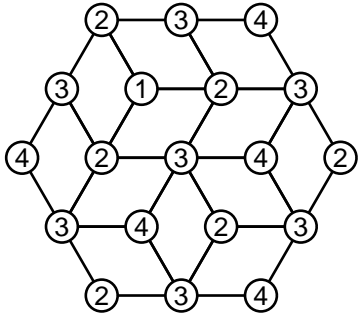}
\end{center}
\caption{The height function corresponding to a tiling.}
\label{fig-heightfun}
\end{figure}

An alternative way to phrase \thmref{mainthm} is in terms of height
functions, which were introduced by Thurston in \cite{Th} as a geometrical
tool for understanding tilings of regions by lozenges or dominos.
A height function is a numerical representation of an individual tiling
of a specified region.  In the case of lozenge tilings of a hexagon,
the vertices of the lozenges occur at points of a certain triangular
lattice that is independent of the particular tiling chosen, and the
height function simply associates a certain integer to each such vertex
so as to describe the shape of the plane partition that corresponds to
the tiling.
Given any lozenge tiling, one
can assign heights $h(u)$ to the points $u$ of the triangular lattice
as follows.  Give the leftmost vertex of the hexagon height $a+c$.
(We choose this height so that the vertex of the bounding box
farthest from the viewer,
which is usually obscured from view,
will have height $0$.)
Suppose that $u$ and $v$ are adjacent lattice points, such that the
edge connecting them does not bisect a lozenge.  If the edge from $u$
to $v$ points directly to the right, set $h(v)=h(u)+1$, and if it
points up and to the right, or down and to the right, set $h(v)=h(u)-1$.
(If it points left, change $+1$ to $-1$ and vice versa.)
If one traces around each lozenge in the
tiling and follows this rule, then every vertex is assigned a height.  It
is not hard to check that the heights are well-defined,
so there is a unique {\it height function} associated to the tiling.
For an example, see Figure~\ref{fig-heightfun}.  Conversely, every way
of assigning integer heights to the vertices of the
triangular grid
that assigns height $a+c$ to the leftmost vertex of the hexagon and that
satisfies the edge relation must be the height function associated to
some unique tiling.  If one views the tiling as a three-dimensional picture
of the solid Young diagram of a plane partition, then the height function
tells how far above a reference plane (of the form $x+y+z=\mbox{constant}$)
each vertex lies.
The values of the height function on the boundary of the hexagon
are constrained, but
all the values in the interior genuinely depend on the tiling.
(It should be mentioned that height functions for lozenge tilings are
implicit in the work of physicists Bl\"ote and Hilhorst,
in the context of the two-dimensional dimer model on a hexagonal lattice;
see \cite{BH}.)

When dealing with normalized coordinates, it is convenient to use a
{\it normalized height function\/}: if we scale the coordinates by
dividing them by $\scalefactor$ (as described above), then we also
divide the height function values by $\scalefactor$.  Also, we define
the average height function to be the average of the height functions
associated with all possible tilings.  (Of course, it is not a height
function itself.)  We will show that almost every height function
closely approximates the average height function, asymptotically:

\begin{thm}
\label{heightfun}
Fix $\alpha, \beta, \gamma > 0$.  As $a,b,c \rightarrow \infty$ with
$a:b:c \rightarrow \alpha:\beta:\gamma$, the normalized average height
function of the $a,b,c$ hexagon converges uniformly to the function
$\Hf_{\alpha,\beta,\gamma}$ with the appropriate
(piecewise linear)
boundary conditions such that
$$
\frac{\partial \Hf_{\alpha,\beta,\gamma}(x,y)}{\partial x}
= 1- 3\Pl_{\alpha,\beta,\gamma}(x,y).
$$
For fixed $\varepsilon > 0$, the probability that the normalized
height function of a random tiling will differ anywhere by more than
$\varepsilon$ from $\Hf_{\alpha,\beta,\gamma}$ is exponentially small
in $\scalefactor^2$, where $\scalefactor =
(a+b+c)/(\alpha+\beta+\gamma)$.
\end{thm}

It is not hard to deduce \thmref{heightfun} (with the exception of the
claim made in the last sentence) from \thmref{mainthm}; in particular,
the differential equation simply results from considering how the
height changes as one crosses lozenges of each orientation.  We will
give the proof in detail in Section~\ref{sec-analysis}.

Unfortunately, although one can recover an explicit formula for
$\Hf_{\alpha,\beta,\gamma}(x,y)$ from the boundary values and the
knowledge of ${\partial \Hf_{\alpha,\beta,\gamma}(x,y)}/{\partial x}$,
we cannot find any simple formula for it.  By contrast, Proposition~17
of \cite{CEP} gives a comparatively simple asymptotic formula for the
average height function for domino tilings of an Aztec diamond.

Our methods also apply to the case of random domino tilings of Aztec
diamonds.  Formula~(7) of Section~4 of \cite{EKLP} is analogous to our
\thmref{countonline}, and can be used in the same way.  It turns out
that the functional arrived at by applying the methods of
Section~\ref{sec-setup} of this paper to that formula is very similar
to the one we will find in Section~\ref{sec-setup}.  After a simple
change of variables, one ends up with essentially the same functional,
but maximized over a slightly different class of functions.  The
methods of Section~\ref{sec-analysis} apply almost without change, and
the methods of Section~\ref{sec-typical} can be adapted to prove
Proposition~17 of \cite{CEP}.  This proof is shorter than the one
given in \cite{CEP}; however, in \cite{CEP}, Proposition~17 comes as a
corollary of a much stronger result (Theorem~1), which the methods of
this paper do not prove.

\section{The Product Formula}
\label{sec-prodform}

In this section, we will prove a refinement of MacMahon's formula, following
methods first used by Elkies et al.\ in \cite{EKLP}.  This refinement
(Proposition~\ref{prodformula}) is not strictly speaking new, since it is
really nothing more than the Weyl dimension formula for finite-dimensional
representations of $\SL(n)$ (we say a few words more about this connection
below).  However, we give our own proof of this result for two reasons: first,
to make this part of the proof self-contained; and second, to illustrate an
expeditious method of proof that has found applications elsewhere (see
\cite{JP} for related formulas derived by the same method).

\begin{figure}
\begin{center}
\leavevmode
\epsfbox[0 0 180 130]{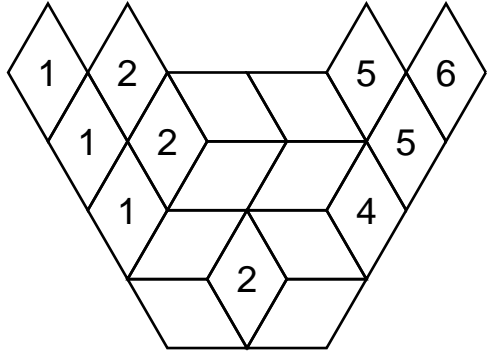}
\end{center}
\caption{The semi-strict Gelfand pattern corresponding to a tiling.}
\label{fig-gelfand}
\end{figure}

Proposition~\ref{prodformula} is stated in terms of Gelfand patterns,
so we must first explain what Gelfand patterns are and what they have
to do with lozenge tilings.  It is not hard to see that a lozenge
tiling of a hexagon is determined by the locations of its vertical
lozenges.  A semi-strict Gelfand pattern is a way to keep track of
these locations.  Specifically, one augments the $a,b,c$ hexagon by
adding $a(a+1)/2$ vertical lozenges on the left and $c(c+1)/2$
vertical lozenges on the right, forming an approximate trapezoid of
upper base $a+b+c$ and lower base $b$, with some triangular
protrusions along its upper border, as in Figure~\ref{fig-gelfand}.
One then associates with each vertical lozenge in the tiling the horizontal
distance from its right corner to the left border of the trapezoid,
which we call its {\it trapezoidal position}.  (When we want to use
the left boundary of the hexagon instead of the left boundary of the
approximate trapezoid as our bench mark, we will speak of the {\it
hexagonal position} of a vertical lozenge.)  For example, consider the
tiling shown in Figure~\ref{fig-shadedhex}; we augment it by adding
$6$ vertical lozenges to form the tiling shown in
Figure~\ref{fig-gelfand}.  The vertical lozenges form rows, and the
only restriction on their placement is that given any two adjacent
vertical lozenges in the same row, there must be exactly one vertical
lozenge between them in the row immediately beneath them.  If we index
the locations of the vertical lozenges with their trapezoidal
positions (the numbers shown in Figure~\ref{fig-gelfand}), we arrive
at the following semi-strict Gelfand pattern:
$$
\begin{array}{ccccccc}
1 &   & 2 &   & 5 &   & 6 \\
&&&&&&\\
  & 1 &   & 2 &   & 5 &   \\
&&&&&&\\
  &   & 1 &   & 4 &   &   \\
&&&&&&\\
  &   &   & 2 &   &   &   \\
\end{array}
$$

In general, a {\it semi-strict Gelfand pattern} is a triangular array of
integers (such as this one), with the property that the entries increase
weakly along diagonals running down and to the right, and the entries increase
strictly along diagonals running up and to the right.  As discussed above,
there is a simple bijection between semi-strict Gelfand patterns with top row
$1,2,\dots,a,a+b+1,a+b+2,\dots,a+b+c$ and lozenge tilings of an equi-angular
hexagon with side lengths $a,b,c,a,b,c$.

Moreover, consider the $k$-th horizontal line from the top in an $a,b,c$
hexagon, where the top edge of the hexagon corresponds to $k=0$, so that the
trapezoidal positions of the vertical lozenges on the $k$-th line are
precisely the entries in the $(k+1)$-st row of the associated semi-strict
Gelfand pattern.  If we discard all lozenges that lie above the line (but
retain all vertical lozenges that are bisected by it and all lozenges that lie
below it), then we get a tiling of a smaller approximate trapezoid, whose
upper border as before consists of triangular protrusions alternating with
straight edges, except that now the protrusions need not be concentrated at
the left and right portions of the border.  The trapezoidal positions of the
vertical lozenges in this tiling are given by the entries in the semi-strict
Gelfand pattern obtained from the original semi-strict Gelfand pattern by
deleting the first $k$ rows.  In fact, if we limit ourselves to tilings of the
$a,b,c$ hexagon in which the locations of the vertical lozenges that are
bisected by the $k$-th horizontal line are specified, then each individual
tiling of this kind corresponds to a {\it pair} of semi-strict Gelfand
patterns.  We have already described one of these Gelfand patterns, which
gives the behavior of the tiling below the cutting line; the other,
which describes the tilings above the line, comes from
reflecting the hexagon through the horizontal axis
(and of course adjoining additional lozenges, as above).
If $0 \le k < \min(a,c)$, then one of the Gelfand
patterns will include some of the augmenting vertical lozenges described
above on both sides (a row of length $a-k$ on one side of the top of
the pattern, and a row of length $c-k$ on the other)
and the other pattern will contain augmenting vertical
lozenges on neither side; if $\min(a,c) \leq k \leq (a+c)/2$, then one
of the Gelfand patterns will contain a row of $|a-k|$ augmenting
lozenges on one side and the other Gelfand pattern will contain
$|c-k|$ on the other side.
(The case $k > (a+c)/2$ is symmetric to the case $k < (a+c)/2$, so
we do not treat it explicitly.)

In \thmref{countonline}, we will use the following formula to
determine how many tilings have a specified distribution of vertical lozenges
along a horizontal line.

\begin{prop}
\label{prodformula}
There are exactly
$$
\prod_{1 \le i < j \le n}\frac{x_j-x_i}{j-i}
$$
semi-strict Gelfand patterns with top row consisting of integers
$x_1,x_2,\dots,x_n$ such that $x_1 < x_2 < \dots < x_n$.
\end{prop}

\propref{prodformula} has an explanation in terms of representation theory.
Gelfand and Tsetlin \cite{GT} show that
semi-strict Gelfand patterns form bases of representations of
$\SL(n)$, and one can deduce \propref{prodformula} from this fact
using the Weyl dimension formula.  (The Gelfand patterns in \cite{GT}
are not semi-strict, but there is an easy transformation that makes
them so: Replace $m_{p,q}$ in equation~(3) of \cite{GT} with
$m_{p,q}+(q-p)$ and then reflect the triangle through its vertical
axis.)  Thus, the novelty of our approach is not that one can count
semi-strict Gelfand patterns, but rather that one can count tilings
with prescribed behavior on a horizontal line (as in
\thmref{countonline}).  Another proof of \propref{prodformula}, and
one that bypasses its representation-theoretic significance, is the
article of Carlitz and Stanley \cite{CS}.  (That article does not deal
directly with semi-strict Gelfand patterns, but it is easy to deduce
\propref{prodformula} from the theorem proved there.)

\begin{proof}
Let $V(x_1,\dots,x_n)$ be the number of semi-strict Gelfand patterns with top
row $x_1,\dots,x_n$.  Given any such pattern, the second row must be of the
form $y_1,\dots,y_{n-1}$ with $x_i \le y_i < x_{i+1}$ for all $i$.  Therefore,
\begin{equation}
\label{Wrec}
V(x_1,\dots,x_n)=\sumL_{y_1=x_1}^{x_2}\sumL_{y_2=x_2}^{x_3}
\dots\sumL_{y_{n-1}=x_{n-1}}^{x_n} V(y_1,\dots,y_{n-1}),
\end{equation}
where the modified summation operator $\sumL$ is defined by
$$
\sumL_{i=s}^t f(i) = \sum_{i=s}^{t-1}f(i).
$$
The advantage to writing it this way is that
\begin{equation}
\label{additive}
\sumL_{i=r}^s f(i) + \sumL_{i=s}^t f(i) = \sumL_{i=r}^t f(i)
\end{equation}
whenever $r < s < t$.  There is a unique way to extend the
definition of
$$
\sumL_{i=s}^t f(i)
$$
to the case $s>t$,
if we require that \eqref{additive}
be true for all $r,s,t$.  Then starting from the base relation
$V(x_1)=1$, we can use \eqref{Wrec} to define $V(x_1,\dots,x_n)$ for
arbitrary integers $x_1,\dots,x_n$ (not necessarily satisfying $x_1 <
\dots < x_n$).

We will prove the formula for $V(x_1,\dots,x_n)$ by induction on $n$.
It is clearly true for $n=1$.  Suppose that for all
$y_1,\dots,y_{n-1}$,
\begin{equation}
\label{indform}
V(y_1,\dots,y_{n-1}) = \prod_{1 \le i < j \le n-1}\frac{y_j-y_i}{j-i}.
\end{equation}
Then $V(y_1,\dots,y_{n-1})$ is a polynomial of total degree
$(n-1)(n-2)/2$ in $y_1,\dots,y_{n-1}$.  When we put \eqref{indform}
into \eqref{Wrec}, we find that after each of the $n$ summations in
\eqref{Wrec}, the result is still a polynomial, and the degree
increases by $1$ each time.  It is easy to check from \eqref{indform}
that the coefficient
of $y_2^1y_3^2\dots y_{n-1}^{n-2}$ in $V(y_1,\dots,y_{n-1})$
is
$$
\frac{1}{1!2!\dots(n-2)!}.
$$
{}From this and \eqref{Wrec}
it follows that the coefficient of $x_2^1x_3^2\dots
x_n^{n-1}$ in $V(x_1,\dots,x_n)$ is
$$
\frac{1}{1!2!\dots(n-1)!}.
$$

We have therefore shown that $V(x_1,\dots,x_n)$ and
\begin{equation}
\label{product}
\prod_{1 \le i < j \le n} \frac{x_j-x_i}{j-i}
\end{equation}
are polynomials in $x_1,\dots,x_n$ of the same total degree, and with
the same coefficient of $x_2^1x_3^2\dots x_n^{n-1}$.  If we can show
that $V(x_1,\dots,x_n)$ is divisible by \eqref{product}, then they
must be equal.  Equivalently, we just need to show that
$V(x_1,\dots,x_n)$ is skew-symmetric in $x_1,\dots,x_n$.

For convenience, let $(s_1,t_1;s_2,t_2;\dots;s_{n-1},t_{n-1})$
denote
$$
\sumL_{y_1=s_1}^{t_1} \sumL_{y_2=s_2}^{t_2} \dots
\sumL_{y_{n-1}=s_{n-1}}^{t_{n-1}} V(y_1,y_2,\dots,y_{n-1}).
$$
We want to show that $(x_1,x_2;x_2,x_3;\dots;x_{n-1},x_n)$ is
skew-symmetric in $x_1,\dots,x_n$.

To start off, notice that for all $i$,
\begin{equation}
\label{summeq}
(\dots;s_i,u_i;\dots) = (\dots;s_i,t_i;\dots) +
(\dots;t_i,u_i;\dots).
\end{equation}
Also, since $V(y_1,\dots,y_{n-1})$ is a skew-symmetric function of
$y_1,\dots,y_{n-1}$ by \eqref{indform}, for $i < j$ we must have
$$
(\dots;s_i,t_i;\dots;s_j,t_j;\dots) =
-(\dots;s_j,t_j;\dots;s_i,t_i;\dots).
$$
The relation $(\dots;s_i,t_i;\dots)=-(\dots;t_i,s_i;\dots)$ follows
easily from the definition of $(\dots;s_i,t_i;\dots)$.
{}From the last two relations, we see that
$(\dots;s_i,t_i;\dots;s_j,t_j;\dots)$ vanishes if
$s_i=s_j$ and $t_i=t_j$, or if $s_i=t_j$ and $s_j=t_i$.

To verify that $(x_1,x_2;\dots;x_{n-1},x_n)$ is skew-symmetric in
$x_1,\dots,x_n$, it suffices to check that it changes sign under
transpositions of adjacent $x_i$'s.  We check the effect of exchanging
$x_{i+1}$ with $x_{i+2}$ as follows.  If we write $x'_k = x_{i+k}$ (to
simplify the subscripts), we have
\begin{eqnarray*}
(\dots;x'_0,x'_2;x'_2,x'_1;x'_1,x'_3;\dots)  & = &
\phantom{\phantom{} + \phantom{}}
(\dots;x'_0,x'_1;x'_2,x'_1;x'_1,x'_2;\dots) \\
& & \phantom{} +
(\dots;x'_0,x'_1;x'_2,x'_1;x'_2,x'_3;\dots) \\
& & \phantom{} +  (\dots;x'_1,x'_2;x'_2,x'_1;x'_1,x'_2;\dots) \\
& & \phantom{} +  (\dots;x'_1,x'_2;x'_2,x'_1;x'_2,x'_3;\dots) .
\end{eqnarray*}
by several applications of \eqref{summeq}.  All terms on the right
except the second are $0$, so
\begin{eqnarray*}
(\dots;x'_0,x'_2;x'_2,x'_1;x'_1,x'_3;\dots)
& = & (\dots;x'_0,x'_1;x'_2,x'_1;x'_2,x'_3;\dots)\\
& = & -(\dots;x'_0,x'_1;x'_1,x'_2;x'_2,x'_3;\dots) .
\end{eqnarray*}
Thus, exchanging $x_{i+1}$ with $x_{i+2}$ introduces a minus sign
whenever $i \ge 1$ and $i+3\le n$.  The other cases (exchanging $x_1$
with $x_2$ or $x_{n-1}$ with $x_n$) are easily dealt with in a similar
way.  Therefore, $V(x_1,\dots,x_n)$ is skew-symmetric in
$x_1,\dots,x_n$, so as discussed above we must have
$$
V(x_1,\dots,x_n) = \prod_{1 \le i < j \le n} \frac{x_j-x_i}{j-i},
$$
as desired.
\end{proof}

Notice that after some manipulation of the product
$$
V(1,2,\dots,a,a+b+1,a+b+2,\dots,a+b+c),
$$
\propref{prodformula} implies MacMahon's formula.  However, our main
application will be to counting tilings with prescribed behavior on horizontal
lines.

Consider the $k$-th horizontal line from the top in an $a,b,c$ hexagon.  If $k
< \min(a,c)$, then in every tiling there must be $k$ vertical lozenges on the
$k$-th line; if $\min(a,c) \le k \le (a+c)/2$, then there must be $\min(a,c)$
vertical lozenges on it.  (As mentioned earlier, symmetry frees us from
needing to treat the case $k > (a+c)/2$, so we will routinely assume $k \leq
(a+c)/2$.)  In either case, note that the number of vertical lozenges on the
$k$-th row is $\min(k,a,c)$.

Define the function
$$
V(x_1,\dots,x_n) = \prod_{1 \le i < j \le n} \frac{x_j-x_i}{j-i},
$$
as in the proof of \propref{prodformula}.  Then we have the following
formulas:

\begin{thm}
\label{countonline}
Suppose we require that the vertical lozenges bisected by
the $k$-th horizontal line from the top in an $a,b,c$ hexagon
occur at hexagonal positions
$1 \le a_1 < a_2 < \dots < a_\ell \le b+\min(k,a)$
(and nowhere else), with $\ell = \min(k,a,c)$.
If $k < \min(a,c)$, there are
$$
V(a_1,a_2,\dots,a_\ell)
V(1,2,\dots,a-k,a-k+a_1,\dots,a-k+a_\ell,
a+b+1,\dots,a+b+c-k)
$$
such tilings.
If $a \leq c$ and $a \le k \le (a+c)/2$, there are
$$
V(1,2,\dots,k-a,k-a+a_1,\dots,k-a+a_\ell)
V(a_1,a_2,\dots,a_\ell,a+b+1,\dots,a+b+c-k)
$$
such tilings
(and a similar formula applies if $c < a$ and $c \le k \le (a+c)/2$).
\end{thm}

For the proof, simply notice that tilings of the parts of the hexagon above
and below the $k$-th line correspond naturally to semi-strict Gelfand patterns
with certain top rows, and then apply \propref{prodformula} to count them.  In
both cases, the two factors correspond to the parts of the tiling that lie
above and below the cutting line.

\section{Setting up the Functional}
\label{sec-setup}

We now turn to the proofs of our main theorems.
As is usually done in situations such as ours,
where one seeks to establish a law of large numbers
for some class of combinatorial objects,
we approach the problem by trying to find the individually most
likely outcome (in this case, the individually most likely behavior
of the height function on a fixed horizontal line), and showing
that outcomes that differ substantially from it are very unlikely,
even considered in aggregate.  We will begin in this section by
setting up a functional to be maximized;  the function that maximizes
it will be a simple transformation of the average height function.

Our method is to focus on the locations of the vertical lozenges rather than
the height function per se.  The two are intimately related, because, as
we move across the tiling from left to right, the (unnormalized) height
decreases by 2 when we cross a vertical lozenge and increases by 1 when we
fail to cross a vertical lozenge.
Thus, in determining the likely locations of vertical lozenges, we
will in effect be solving for the average height function.
\thmref{countonline} tells us that we can count tilings with
prescribed behavior on horizontal lines, so we will start off
by taking the logarithm of the product formula in \thmref{countonline}
and interpreting it as a Riemann sum for a double integral.

In fact, it will be convenient to look first at a more
general product, and then apply our analysis of it to the product
appearing in \thmref{countonline}.
Fix positive integers $\ell$ and $\nmiddle$ satisfying $\ell \le
\nmiddle$, and non-negative integers $\nleft$ and $\nright$.  We will
try to determine the distribution of $(a_1,a_2,\dots,a_\ell)$
satisfying $1 \le a_1 < a_2 < \dots < a_\ell \le \nmiddle$ that
maximizes
$$
V(1,2,\dots,\nleft,\nleft+a_1,\nleft+a_2,\dots,
\nleft+a_\ell,\nleft+\nmiddle+1,\nleft+\nmiddle+2,\dots,
\nleft+\nmiddle+\nright).
$$
For convenience, let $b_i$ denote the $i$-th element of the sequence
$1,2,\dots,\nleft,\nleft+a_1,\dots,\nleft+a_\ell,\nleft+\nmiddle+1,\dots,
\nleft+\nmiddle+\nright$.

Set $\rhL = \nleft/\nmiddle$, $\lambda = \ell/\nmiddle$, and
$\rhR = \nright/\nmiddle$.
We will work in the limit as $\nmiddle \rightarrow \infty$,
with $\rhL$, $\lambda$, and $\rhR$
tending toward definite limits.
For a more convenient way to keep track of $a_1,\dots,a_\ell$ as we
pass to the continuum limit, we define a
weakly increasing function $A : [0,1] \rightarrow [0,1]$ as follows.
Set $A(i/\nmiddle) = \#\{j : a_j \le i\}/\nmiddle$ for $0 \le i \le \nmiddle$,
and then interpolate by straight lines between these points.
Similarly, set $B(i/\nbig) = \#\{j : b_j \le i\}/\nbig$ for $0 \le i \le
\nbig,$
where $\nbig \stackrel{\textup{def}}{=} \nleft+\nmiddle+\nright$, and then
interpolate by straight lines.
To simplify the notation later, set $\rhoLR = \rhL + \rhR$,
so that
$\nbig = (1+\rhoLR)\nmiddle$,
$A(0) = 0$,
$A(1) = \lambda$,
$B(0) = 0$, and
$B(1) = (\lambda+\rhoLR)/(1+\rhoLR)$.

The functions $A$ and $B$ satisfy the Lipschitz condition with constant
1; that is, $|A(s)-A(t)|$ and $|B(s)-B(t)|$ are bounded by $|s-t|$.
Note that the derivatives $A'$ and $B'$ are not defined everywhere,
but they are undefined only at isolated points, and where they are
defined they equal either $0$ or $1$; when we make statements about
$A'$ and $B'$, we will typically ignore the points of
non-differentiability.  We can also derive a simple relation between
$A'$ and $B'$.  To do so, notice that for $0 \le i \le \nmiddle$, we
have
$$
B\left(\frac{i+\nleft}{\nbig}\right)
= \frac{\#\{j : b_j \le i+\nleft \}}{\nbig}
= \frac{\#\{j : a_j \le i\}}{\nbig}+\frac{\nleft}{\nbig}
=A\left(\frac{i}{\nbig}\cdot\frac{\nbig}{\nmiddle}\right)
\frac{\nmiddle}{\nbig} + \frac{\nleft}{\nbig};
$$
if we set $t=i/\nbig$, we find that this equation becomes
\begin{equation}
\label{BArel}
B\left(t+\frac{\rhL}{1+\rhoLR}\right) =
\frac{A\left((1+\rhoLR) t\right)}{1+\rhoLR}+\frac{\rhL}{1+\rhoLR}.
\end{equation}
For $0 \le t \le
1/(1+\rhoLR)$, the values of $A$ and $B$ occurring in
\eqref{BArel} are defined by interpolating linearly between points at
which we have just shown that \eqref{BArel} holds, so it must hold for
all such $t$.  Therefore, for $0 \le t \le 1/(1+\rhoLR)$, we have
$$
B'\left(t+\frac{\rhL}{1+\rhoLR}\right) = A'((1+\rhoLR)t),
$$
except at isolated points of non-differentiability.  All other values
of $B'$ are $1$, since it follows immediately from the definition of
$B$ that for $0 < t < \rhL/(1+\rhoLR)$ or $1 > t >
(1+\rhL)/(1+\rhoLR)$, we have $B'(t)=1$.

We have $a_i = \nmiddle A^{-1}(i/\nmiddle)$, $b_i = \nbig
B^{-1}(i/\nbig)$.  (Whenever we refer to $A^{-1}(t)$, we consider
it to take the smallest value possible, to avoid ambiguity, and we
interpret $B^{-1}(t)$ analogously.)  When we take the logarithm and
then multiply by $\nmiddle^{-2}$, the double product
$$
V(b_1,\dots,b_{\nleft+\ell+\nright}) =
\prod_{1 \leq i < j \leq \nleft+\ell+\nright} \frac{b_j-b_i}{j-i}
$$
ought to approach an integral like
$$
\left(1+\rhoLR\right)^2 \iint_{0 \leq s < t \leq
\frac{\lambda+\rhoLR}{1+\rhoLR}} \log \frac{B^{-1}(t)-B^{-1}(s)}{t-s} \,ds
\,dt.
$$
(The factor of $(1+\rhoLR)^2$ appears because
$\nbig=(1+\rhoLR)\nmiddle$ and we rescaled by $\nmiddle^{-2}$ instead
of $\nbig^{-2}$.)  In the appendix, we will justify this claim
rigorously, except with the function $B$ replaced by a nicer function
$C$.  (The justification is not very difficult, but it is long enough
that here it would be a distraction.)

The conclusion from the appendix is that
$$
\frac{\log V(b_1,\dots,b_{\nleft+\ell+\nright})}{\nmiddle^2}
=
(1+\rhoLR)^2 \iint_{0 \leq s < t \leq \frac{\lambda+\rhoLR}{1+\rhoLR}}
\log \frac{C^{-1}(t)-C^{-1}(s)}{t-s} \,ds \,dt + o(1),
$$
where $C$ is a certain strictly increasing, continuous, piecewise
linear function satisfying $C(0) = 0$, $C(1) =
(\lambda+\rhoLR)/(1+\rhoLR)$, and $|B'-C'| = O(1/\nbig)$.  Note that
$C$ has a continuous inverse (unlike $B$, which is only weakly
increasing and thus may not even have an inverse).

By symmetry, the integral equals
$$
\frac{(1+\rhoLR)^2}{2}
\int_0^{\frac{\lambda+\rhoLR}{1+\rhoLR}}
\int_0^{\frac{\lambda+\rhoLR}{1+\rhoLR}}
\log \frac{C^{-1}(t)-C^{-1}(s)}{t-s} \,ds \,dt.
$$
We can write the integral as
\begin{eqnarray*}
& & \frac{(1+\rhoLR)^2}{2} \int_0^{\frac{\lambda+\rhoLR}{1+\rhoLR}}
\int_0^{\frac{\lambda+\rhoLR}{1+\rhoLR}}
\log \frac{|C^{-1}(t)-C^{-1}(s)|}{|t-s|} \,ds \,dt \\
& = & \frac{(1+\rhoLR)^2}{2} \int_0^{\frac{\lambda+\rhoLR}{1+\rhoLR}}
\int_0^{\frac{\lambda+\rhoLR}{1+\rhoLR}}
\log |C^{-1}(t)-C^{-1}(s)| \,ds \,dt \\
& & \phantom{}
- \frac{(1+\rhoLR)^2}{2} \int_0^{\frac{\lambda+\rhoLR}{1+\rhoLR}}
\int_0^{\frac{\lambda+\rhoLR}{1+\rhoLR}}
\log |t-s| \,ds \,dt ;
\end{eqnarray*}
since the integral being subtracted is a constant, we can ignore it.
(The individual integrands are unbounded, but since the singularities
are merely logarithmic they do not interfere with integrability.)
Letting $u=C^{-1}(s)$ and $v=C^{-1}(t)$, we can rewrite the part
that matters as
$$
\frac{(1+\rhoLR)^2}{2}
\int_0^1 \int_0^1 C'(u) C'(v) \log |u-v| \,du \,dv.
$$
Since $|C'-B'| = O(1/\nbig)$, this integral differs by
$o(1)$ from
\begin{equation}
\label{inset}
\frac{(1+\rhoLR)^2}{2}
\int_0^1 \int_0^1 B'(u) B'(v) \log |u-v| \,du \,dv.
\end{equation}

We will now use our formula expressing $B'$ in terms of $A'$.
Recall that
$$
B'\left(\frac{t+\rhL}{1+\rhoLR}\right) = A'(t)
$$
for $0 \le t \le 1$, and that $B'(t)=1$ for $t \in (-\rhL,0) \cup
(1,1+\rhR)$.  To take advantage of this, we change variables to $s$
and $t$ (which are different from the $s$ and $t$ used earlier in the
article) with $u = (s+\rhL)/(1+\rhoLR)$ and $v = (t+\rhL)/(1+\rhoLR)$.
Then when $0 \le s,t \le 1$ we have $B'(u)=A'(s)$ and $B'(v)=A'(t)$;
elsewhere $B'$ is $1$.  Thus, \eqref{inset} is equal to
\begin{equation}
\label{almostdone}
\frac{1}{2} \int_{-\rhL}^{1+\rhR}
\int_{-\rhL}^{1+\rhR}
(A'(s)+I(s))(A'(t)+I(t)) \log \left|\frac{s-t}{1+\rhoLR}\right| \,ds
\,dt,
\end{equation}
where $I$ is the characteristic function of $[-\rhL,0] \cup
[1,1+\rhR]$, and where we set $A'=0$ outside $[0,1]$.

Removing the $1+\rhoLR$ from the denominator of the argument of the
logarithm simply adds
$$
\frac{1}{2}\int_{-\rhL}^{1+\rhR}
\int_{-\rhL}^{1+\rhR}
(A'(s)+I(s))(A'(t)+I(t)) \log |1+\rhoLR| \,ds \,dt
$$
to the integral;  this quantity is a constant because
the only occurrence of $A$ in it is through the integral
$$
\int_{-\rhL}^{1+\rhR} A'(s) \,ds = A(1)-A(0) = \lambda
$$
(and the square of this integral).
We can also change the range of integration in \eqref{almostdone} to
the entire plane
(since the integrand has support only
in the rectangle $[-\rhL,1+\rhR] \times [-\rhL,1+\rhR]$).
Thus, we have arrived at the result that, for some
irrelevant constant $K$,
$\nmiddle^{-2}\log V(b_1,\dots,b_{\nleft+\ell+\nright})$ equals
\begin{equation}
\label{functional}
\frac{1}{2} \intr \intr (A'(s)+I(s))(A'(t)+I(t))
\log |s-t| \,ds \,dt + K + o(1).
\end{equation}

We can now apply this to \thmref{countonline}.  Suppose that on the
$k$-th line from the top in the $a,b,c$ hexagon (with $k \le
(a+c)/2$), the vertical lozenges have hexagonal positions $1 \le a_1 <
a_2 < \dots < a_\ell \le \nmiddle$, where $\nmiddle
\stackrel{\textup{def}}{=} b+\min(k,a,c)$.  Define the function $A$ as
above.  Then our analysis so far, combined with \thmref{countonline},
shows that if we take the logarithm of the number of tilings with the
given behavior on the $k$-th line, and divide by $\nmiddle^2$, then we
get the sum of two terms of the form \eqref{functional}; for,
\thmref{countonline} gives us a product of two $V$-expressions (whose
exact nature depends on whether $k$ lies between $a$ and $c$), and
when we take logarithms and divide by $\nmiddle^2$, we get two terms,
each of which is half of a double integral (plus negligibly small
terms and irrelevant constants).

To put this into an appropriately general context, define the bilinear
form $\langle,\rangle$ by
\begin{equation}
\label{inproddef}
\inprod{f}{g} = \intr \intr f'(x) g'(y) \log |x-y| \,dx \,dy
\end{equation}
for suitable functions $f$ and $g$ (for our purposes, functions such
that their derivatives exist almost everywhere, are bounded, are
integrable, and have compact support).  We will now use this notation
to continue the analysis begun in the previous paragraph.

In this paragraph we will systematically omit additive constants and
$o(1)$ terms, since they would be a distraction.  If $k < \min(a,c)$,
then one of the two terms derived from \thmref{countonline} is
$\frac{1}{2}\inprod{A}{A}$, and the other is
$\frac{1}{2}\inprod{A+J_0}{A+J_0}$, where $J_0$ is a continuous
function with derivative equal to the characteristic function of
$[-|a-k|/\nmiddle,0] \cup [1,1+|c-k|/\nmiddle]$.  If $\min(a,c) \le k
\le (a+c)/2$, then one term is $\frac{1}{2}\inprod{A+J_1}{A+J_1}$ and
the other is $\frac{1}{2}\inprod{A+J_2}{A+J_2}$, where the derivative
of $J_1$ is the characteristic function of $[-|a-k|/\nmiddle,0]$ and
that of $J_2$ is the characteristic function of
$[1,1+|c-k|/\nmiddle]$.

Now a few simple algebraic manipulations bring these results to the
following form.

\begin{prop}
\label{functionalprop}
Let $\nmiddle = b + \min(k,a,c)$, $\rhL = |a-k|/\nmiddle$, and $\rhR =
|c-k|/\nmiddle$.  Then the logarithm of the number of tilings of an
$a,b,c$ hexagon with vertical lozenges at hexagonal positions
$a_1,\dots,a_\ell$ (and nowhere else) along the $k$-th line (where
$\ell = \min(k,a,c)$ and $1 \le a_1 < \cdots < a_\ell \le \nmiddle$),
when divided by $\nmiddle^2$, equals
$$
\inprod{A+J}{A+J}+\textup{constant}+o(1),
$$
where $J$ is any continuous function whose derivative is half the
characteristic function of $[-\rhL,0] \cup [1,1+\rhR]$, and $A$ is
defined (as earlier) by interpolating linearly between the values
$A(i/n) = \#\{j : a_j \le i\}/n$ for $0 \le i \le n$.
\end{prop}

Recall that we are working in the limit as $n \rightarrow \infty$ with
$\rhL$, $\rhR$, and $\lambda$ converging towards fixed limits.  It is easy
to check that the unspecified constant in
Proposition~\ref{functionalprop} converges as well.
We would like everything to be stated in terms of the limiting values
of $\rhL$, $\rhR$, and $\lambda$.  Replacing $\rhL$ and $\rhR$ in the
definition of $J$ by their limits changes $\inprod{A+J}{A+J}$ by only
$o(1)$, and we can increase $A'$ by $o(1)$ in such a way that $A(0)$
remains $0$, $A(1)$ becomes the limiting value of $\lambda$, and $A$
still satisfies the Lipschitz condition.  Thus, we can let $\rhL$,
$\rhR$, and $\lambda$ denote their limiting values from now on.

We have now re-framed our problem.  We must find a function $A$ that
maximizes $\V(A) \stackrel{\textup{def}}{=} \inprod{A+J}{A+J}$,
subject to certain conditions.  We will look at real-valued functions
$A$ on $[0,1]$ that are continuous, weakly increasing, and subject to
the following constraints: $A(0)=0$, $A(1)=\lambda$, and $A$ must
satisfy a Lipschitz condition with constant $1$ (so $0 \le A' \le 1$
where $A'$ is defined).  For convenience, define $A(t) = 0$ for $t <
0$ and $A(t) = \lambda$ for $t > 1$.  Call a function $A$ that
satisfies these conditions {\it admissible}. Clearly, the functions
$A$ considered in this section are admissible.  We will show in the
next section that there is a unique admissible function $A$ that
maximizes $\V(A)$.  (Notice that every admissible function $A$ is
differentiable almost everywhere, and $A'$ is integrable and has
compact support; for a proof of the necessary facts from real
analysis, see Theorem~7.18 of \cite{R}.  Thus, $\V(A)$ makes sense for
every admissible $A$.)

\section{Analyzing the Functional}
\label{sec-analysis}

Let $\Admis$ be the set of admissible functions.  We can topologize
$\Admis$ using the sup norm, $L^1$ norm, or $L^2$ norm on $[0,1]$; it
is straightforward to show that they all give the same topology, and that
$\Admis$ is compact.  In this section, we will show that $\V$ is a
continuous function on $\Admis$, so it must attain a maximum.  We will
show furthermore that there is a unique function $A \in \Admis$ such
that $\V(A)$ is maximal.

The proof will use several useful formulas for the bilinear form
$\langle, \rangle$ (formulas \eqref{innerprod1} to \eqref{negdef}).
These formulas are derived in \cite{LS}; we repeat the derivations
here for completeness.  One can find similar analysis in \cite{VK1}
and \cite{VK2}.

The formulas are stated in terms of the Fourier and Hilbert
transforms.  For sufficiently well-behaved functions $f$, define the
Fourier transform $\widehat{f}$ of $f$ by
$$
\widehat{f}(x) = \intr f(t) e^{-ixt} \,dt
$$
and the Hilbert transform $\widetilde{f}$ by
$$
\widetilde{f}(x) = \frac{1}{\pi} \intr \frac{f(t)}{x-t} \,dt
$$
(which we make sense of by taking the Cauchy principal value).  Note
that the Fourier transform of $\widetilde{f}$ is $x \mapsto i(\sgn
x)\widehat{f}(x)$, and that of $f'$ is $x \mapsto i x \widehat{f}(x)$;
multiplying by $i (\sgn x)$ and multiplying by $i x$ commute with each
other, so differentiation commutes with the Hilbert transform.

Integration by parts with respect to $y$ in \eqref{inproddef} shows
that
\begin{equation}
\label{innerprod1}
\inprod{f}{g} = \pi \intr f'(x) \widetilde{g}(x) \,dx.
\end{equation}
When done with respect to $x$, it also shows that
\begin{equation}
\label{innerprod2}
\inprod{f}{g} = - \pi \intr f(x) \widetilde{g'}(x) \,dx.
\end{equation}
Unfortunately, Hilbert transforms are not always defined.  For our
purposes, it is enough to note that \eqref{innerprod1} makes sense and
is true when $f'$ and $g$ have compact support, and similarly that
\eqref{innerprod2} holds when $f$ and $g'$ have compact support.

If we set $g=f$ and apply Parseval's identity to \eqref{innerprod1},
we find that when $f$ has compact support,
\begin{equation}
\label{negdef}
\inprod{f}{f} = -\frac12 \intr |\widehat{f}(x)|^2 |x| \,dx.
\end{equation}
Thus, the bilinear form $\langle,\rangle$ is negative definite (on
functions of compact support).

We can now prove easily that $\V$ is a continuous function on
$\Admis$.  To do so, notice that the definition of $\V$ and
\eqref{innerprod1} imply that
$$
\V(A_1)-\V(A_2) = \inprod{A_1+A_2+2J}{A_1-A_2} =
\pi \intr f_1'(x) \widetilde{f_2}(x) \,dx,
$$
where $f_1 = A_1+A_2+2J$ and $f_2=A_1-A_2$.  Thus, since $|f_1'(x)|
\le 2$ for all $x$ and $f_1' = 0$ outside some interval $I$ not
depending on $A_1$ and $A_2$,
$$
|\V(A_1)-\V(A_2)| \le 2\pi \int_I |\widetilde{f_2}|
\ll \left( \int_I |\widetilde{f_2}|^2 \right)^{1/2}
\le ||\widetilde{f_2}||_2.
$$
(The second bound follows from applying the Cauchy-Schwartz
inequality to $|\widetilde{f_2}| \cdot 1$.)

It is known (see Theorem~90 of \cite{T}) and easy to prove (combine
Parseval's identity with the formula for the Fourier transform of a
Hilbert transform) that $||\widetilde{f_2}||_2 =
||f_2||_2$.  Thus, $|\V(A_1)-\V(A_2)| =
O(||A_1-A_2||_2)$, so the function $\V$ is continuous on $\Admis$.

Because $\Admis$ is compact, $\V$ must attain a maximum on $\Admis$.
Now we apply the identity
$$
\frac{\V(A_1)+\V(A_2)}{2} = \V\left(\frac{A_1+A_2}{2}\right)
+ \left\langle \frac{A_1-A_2}{2},\frac{A_1-A_2}{2}
\right\rangle,
$$
which is a form of
the polarization identity for quadratic
forms.  Because $(A_1-A_2)/2$ has compact support, \eqref{negdef}
implies that
$$
\frac{\V(A_1)+\V(A_2)}{2} \le \V\left(\frac{A_1+A_2}{2}\right),
$$
with equality if and only if $A_1=A_2$.  Thus, two different
admissible functions could not both maximize $\V$, since then their
average would give an even larger value.  Therefore, there is a unique
admissible function that maximizes $\V$.  Let $\A$ be that function.
(Notice that $\A$ depends on $\lambda$, $\rhL$, and $\rhR$, and hence
on $\alpha$, $\beta$, and $\gamma$, although our notation does not
reflect this dependence.)

We are now almost at the point of being able to prove that there is a
function $\Hf_{\alpha,\beta,\gamma}$ such that \thmref{heightfun} holds
(except for the part relating $\Hf_{\alpha,\beta,\gamma}$ to the explicitly
given function $\Pl_{\alpha,\beta,\gamma}$).
First, we need to relate $A$
to the normalized average height function.

Assume that $k/(a+c) \rightarrow \kappa$ as $a,b,c \rightarrow \infty$
for some $\kappa$ satisfying $0 \leq \kappa \leq 1$.  Choose
normalized coordinates for the $a,b,c$ hexagon so that the $k$-th
horizontal line from the top has normalized length $1$, and in
particular coordinatize that line so that its left endpoint is $0$ and
its right endpoint is $1$; equivalently, coordinatize the
$\alpha,\beta,\gamma$ hexagon so that the horizontal line that cuts it
proportionately $\kappa$ of the way from its upper border has length
1.  (The truth or falsity of \thmref{heightfun} is clearly unaffected
by our choice of coordinates.)  We can then identify the scaling
factor $\scalefactor$ with $\nmiddle$, to within a factor of $1+o(1)$.
Given a tiling of the $a,b,c$ hexagon, if we scan to the right along
this line, whenever we cross a vertical lozenge the normalized height
function decreases by $2/\nmiddle$, and whenever we cross a location
that could hold a vertical lozenge but does not the normalized height
function increases by $1/\nmiddle$.  It follows that the normalized
height function at location $t$ is given by
$$
-2A(t)+(t-A(t)) = t-3A(t)
$$
plus the value at $t=0$, since this function changes by the same amount
as the normalized height function does as one moves to the right.

Let $\varepsilon > 0$.  Then there exists a $\delta > 0$ such that
$|\V(\A)-\V(A)| < \delta$ implies $\sup_t |\A(t)-A(t)| < \varepsilon$.
(This claim holds for every continuous function on a compact space
that takes its maximal value at a unique point.)  For $\nmiddle$
sufficiently large, \propref{functionalprop} implies that if $\sup_t
|\A(t)-A(t)| \ge \varepsilon$, then in a random tiling, every behavior
within $o(1)$ of $\A$ is at least $e^{\nmiddle^2(\delta+o(1))}$ times
as likely to occur along the $k$-th line as the behavior $A$ is.
Since the number of possibilities for $A$ is only exponential in
$\nmiddle$, the probability that $\sup_t |\A(t)-A(t)| \ge \varepsilon$
is exponentially small in $\nmiddle^2$ (and hence in $\sigma^2$).
In other words, the probability that a random height function differs
along this line from the height function obtained from $\A$ is
exponentially small in $\sigma^2$.  It follows that the same is true
without the restriction to the particular horizontal line, because of
the Lipschitz constraint on height functions:  if we show that large
deviations are unlikely on a sufficiently dense (but finite) set of
horizontal lines, then the same holds even between them.
Thus, we have proved
\thmref{heightfun}, except for the connection between
$\Hf_{\alpha,\beta,\gamma}$ and $\Pl_{\alpha,\beta,\gamma}$.

Furthermore, the density of vertical lozenges near location $t$ along
a given horizontal line is almost always approximately equal to
$\A'(t)$.  We can make this claim precise and justify it as follows.
Given a random tiling, $A(t)$ gives the number of vertical lozenges to
the left of the location $t$, divided by $\nmiddle$ (plus
$O(\nmiddle^{-1})$ if $t$ is not at a vertex of the
underlying triangular lattice).  Thus, the number of vertical lozenges
in an interval $[a,b]$ is $\nmiddle(A(b)-A(a)) + O(1)$.  We have seen
that the probability that this quantity will differ by more than
$\varepsilon \nmiddle$ {}from $\nmiddle(\A(b)-\A(a))$ is exponentially
small in $\nmiddle^2$.  Therefore, as $\nmiddle \rightarrow \infty$
(equivalently, $\scalefactor \rightarrow \infty$), the expected value
of $A(b)-A(a)$ is $\A(b)-\A(a)+o(1)$.  Thus, we get the expected
number of vertical lozenges in $[a,b]$, which is also equal to the
expected number of vertical lozenges in $(a,b)$ (up to a negligible
error).  Notice that the $o(1)$ error term is uniform in the choice of
$a$, $b$, and
the horizontal line, because the probability of a large deviation in
height anywhere in the hexagon is small.

If we take the result we have just proved for horizontal line segments
and integrate it over the horizontal line segments that constitute the
interior of any smooth simple closed curve, then we can conclude that
\thmref{mainthm} holds, except for the explicit determination of
$\Pl_{\alpha,\beta,\gamma}$ (which is equivalent to the explicit
determination of $\A$, since $\Pl = \A'$).  Also, notice that our
method of proof implies that $\Hf_{\alpha,\beta,\gamma}$ and
$\Pl_{\alpha,\beta,\gamma}$ must satisfy
$$
\frac{\partial \Hf_{\alpha,\beta,\gamma}(x,y)}{\partial x}
= 1-3\Pl_{\alpha,\beta,\gamma}(x,y),
$$
as desired (although it does not yet determine them explicitly).
Thus, all that remains to be done is to determine the maximizing
function $\A$ explicitly.  We will do so in Section~\ref{sec-typical}.

\section{The Typical Height Function}
\label{sec-typical}

Unfortunately, it is not clear how to find the admissible function $\A$ that
maximizes $\V(\A)$.  Ordinary calculus of variations techniques will not
produce admissible solutions.  However, we will see that techniques similar to
those used in \cite{LS} and \cite{VK1,VK2} can be used to verify that a
function $\A$ maximizes $\V(\A)$, if we can guess $\A$.  (It is not clear
a priori that the techniques will work, but fortunately everything works out
just as one would hope.)

As we saw in Section~\ref{sec-analysis}, for the cases that are needed in the
proof of \thmref{mainthm}, an explicit formula for $\A$ is equivalent to one
for $\Pl_{\alpha,\beta,\gamma}$; since we know already what the answer should
be, guessing it will not present a problem.  In Section~\ref{sec-intro}, we
tried to give some motivation by presenting a slightly simpler description of
$\Pl_{\alpha,\beta,\gamma}$ than the explicit formula.  However, we do not
know of any straightforward way to guess the answer from scratch.  We
arrived at it
partly by analogy with the arctangent formula for random domino tilings of
Aztec diamonds (Theorem~1 of \cite{CEP}), partly on the basis of symmetry and
simplicity, and partly on the basis of numerical evidence.

To avoid unnecessarily complicated notation, we will solve the problem
in greater generality than is needed simply for \thmref{mainthm}.  We
will deal with the case of arbitrary $\lambda$ satisfying $0 < \lambda
< 1$, and arbitrary non-negative $\rhL$ and $\rhR$ (which we assume
for simplicity are strictly greater than $0$).  We will use
the same notation as earlier in the paper; for example, we set $\rhoLR =
\rhL + \rhR$.  Of course, guessing the admissible function that
maximizes the functional in general requires additional effort, but
the symmetry and elegance of the general formulas are a helpful guide.

We will express the maximizing function $\A$ in terms of auxiliary
functions $f_1$ and $f_2$.  Define
$$
f_1(t) = 2t(1-t) - (\lambda^2 + \rhoLR\lambda-\rhR)t -
(\lambda^2+\rhoLR\lambda-\rhL)(1-t)
$$
and
$$
f_2(t) = (\rhoLR+2\lambda)^2t(1-t) - (\lambda^2+\rhoLR\lambda-\rhR)^2t -
(\lambda^2+\rhoLR\lambda-\rhL)^2(1-t).
$$
(Note that both expressions are invariant under
$(t,\rhL,\rhR) \leftrightarrow (1-t,\rhR,\rhL)$;
this observation reduces some of the work involved
in verifying the claims that follow.)
Since the discriminant of $f_2(t)$ is
$$
16\lambda(1-\lambda)(\lambda+\rhR)(\lambda+\rhL)
(\lambda+\rhoLR)(\lambda+\rhoLR+1),
$$
$f_2(t)$ has distinct real roots $r_1 < r_2$.
We can show that both roots are in $[0,1]$ as follows.
Since $f_2(0)$ and $f_2(1)$ are at
most $0$, both roots of $f_2$ lie in $[0,1]$ if the point at which $f_2$
achieves its maximum does.  The maximum occurs when
$$
t = \frac{(1-\lambda)\rhL^2+\lambda\rhR^2 + 2\lambda^2
+(2-\lambda)\lambda\rhL + (2+\lambda)\lambda\rhR + \rhL\rhR}
{(\rhL+\rhR+2\lambda)^2},
$$
and this point is easily seen to lie in $[0,1]$.

We will specify the function $\A$ by specifying its derivative $\A'$,
which together with the condition $\A(0)=0$ uniquely determines the
function.  (We will then check that the newly defined $\A$ maximizes
$\V$, and is thus the same as the previous $\A$.)  For $t \in
(r_1,r_2)$ define
$$
\A'(t) = \frac{1}{\pi}\cot^{-1}
\left(\frac{f_1(t)}{\sqrt{f_2(t)}}\right).
$$
For $t \in [0,r_1]$ define $\A'(t) = \lim_{t \rightarrow r_1+} \A'(t)$, and
similarly for $t \in [r_2,1]$ define $\A'(t) = \lim_{t \rightarrow r_2-}
\A'(t)$.  We can show that the first limit will be $0$ or $1$ if $r_1 \in
(0,1)$, and the second will be $0$ or $1$ if $r_2 \in (0,1)$; to verify this,
it suffices to check (using resultants, for example) that $f_1$ and $f_2$
cannot have a common root in $(0,1)$, from which it follows that at $r_1$ or
$r_2$ the denominator of the argument of the arccotangent vanishes without the
numerator vanishing.  Notice that $\A'(0) = 0$ if $f_1(r_1) > 0$, and
$\A'(0)=1$ if $f_1(r_1) < 0$.  Similarly, $\A'(1) = 0$ if $f_1(r_2) > 0$, and
$\A'(1)=1$ if $f_1(r_2) < 0$.  Also, if $r_1=0$, then $f_1$ and $f_2$ both
vanish at $0$, and it follows that $\lim_{t \rightarrow 0+} \A'(t) = \frac12$;
similarly, if $r_2=1$ then $\lim_{t \rightarrow 1-} \A'(t) = \frac12$.

Let $\A : [0,1] \rightarrow [0,1]$ be the unique function satisfying
$\A(0)=0$ with derivative $\A'$.  We will show later in this section
that $\A(1)=\lambda$, from which it follows that $\A$ is an admissible
function (since the other conditions are clearly satisfied).  (We then
extend $\A$ to a function on all of $\R$ in the usual way, so that
$\A(t) = \lambda$ for $t > 1$ and $\A(t)=0$ for $t < 0.$) We will
prove that $\A$ is the unique admissible function such that $\V(\A)$
is maximal.

Before beginning the proof, it is helpful for motivation
to examine what the calculus of
variations tells us about how the maximizing function should
behave.
For every admissible function $A$ we have
$$
\V(A) = \inprod{A+J}{A+J},
$$
where $J$ is any continuous function whose derivative is half the
characteristic function of $[-\rhL,0] \cup [1,1+\rhR]$.
Suppose we perturb our function $\A$ by adding
to it a small bump centered at $t$, which we write as $\varepsilon
\delta_t$
where $\delta_t$ is a delta function.  (We should actually perturb by
a bump rather than a spike, in order to try to maintain the Lipschitz
constraint as long as $\varepsilon$ is small enough;  however,
treating $\delta_t$ as a delta function gives the right answer.)
Because
$$
\inprod{\A+J+\varepsilon\delta_t}{\A+J+\varepsilon\delta_t}
= \inprod{\A+J}{\A+J} + 2\varepsilon \inprod{\delta_t}{\A+J} +
O(\varepsilon^2),
$$
the first variation of $\V$ is $2\inprod{\delta_t}{\A+J}$.  By
\eqref{innerprod2},
$$
\inprod{\delta_t}{\A+J} = -\pi \intr \delta_t(x) \left(
\widetilde{\A'}(x) + \widetilde{J'}(x)\right) \, dx
= -\pi\left(\widetilde{\A'}(t) + \widetilde{J'}(t)\right).
$$
Thus, in order for the first variation to vanish,
we must have
$\widetilde{\A'}(t) + \widetilde{J'}(t) = 0$.  When $\A'(t) \in (0,1)$
this must be the case, assuming $\A$ maximizes $\V(\A)$.
However, when $\A'(t) \in \{0,1\}$, every perturbation violates
the Lipschitz constraint, and we can conclude nothing.

Our strategy for proving that $\A$ maximizes $\V(\A)$ will take
advantage of the vanishing of
$\widetilde{\A'}(t) + \widetilde{J'}(t)$ (which we will prove
directly later in this section).
To begin, for any admissible $A$ we write
$$
\V(A)-\V(\A) = \inprod{A-\A}{A-\A} + 2\inprod{A-\A}{\A+J}.
$$
Because $\inprod{A-\A}{A-\A} \le 0$ with equality if and only if
$A=\A$,
in order to prove that $\V(A) \le \V(\A)$,
we need only show that
$$
2\inprod{A-\A}{\A+J} \le 0.
$$
To show that this is the case, we start by using
\eqref{innerprod2}, which
tells us that
$$
\inprod{A-\A}{\A+J} =
-\pi \intr (A(t)-\A(t)) (\widetilde{\A'}(t)+
\widetilde{J'}(t)) \,dt.
$$
Thus, we want to show that
\begin{equation}
\label{pos}
\intr (A(t)-\A(t))(\widetilde{\A'}(t)+\widetilde{J'}(t)) \,dt \ge 0.
\end{equation}
We will prove that the integrand is everywhere non-negative, by showing
that $\widetilde{\A'}(t) + \widetilde{J'}(t) = 0$ when $t \in
(r_1,r_2)$ (the interval where $\A'(t) \in (0,1)$), and that in the
rest of $(0,1)$ its sign is the same as that of $A(t)-\A(t)$.
(Outside $(0,1)$, $A(t) = \A(t)$.)

In order to prove these facts, we will apply the following theorem.
For a proof, see Theorem~93 of \cite[p.~125]{T}.

\begin{thm}
\label{titchthm}
Let $\Phi$ be a holomorphic function on the upper half plane, such
that the integrals
$$
\intr |\Phi(x+iy)|^2 \, dx
$$
exist for all $y > 0$ and are bounded.  Then there exists a function
$\Philim$ on the real line such that for almost all $x$, $\Phi(x+iy)
\rightarrow \Philim(x)$ as $y \rightarrow 0+$, and the imaginary part
of $\Philim$ is the Hilbert transform of the real part of $\Philim$.
\end{thm}

We will apply this theorem to determine the Hilbert transform of
$\A'+J'$.    To prepare for
the application of the theorem, we begin by defining, for $t \in
(r_1,r_2)$,
$$
g(t) = \frac{f_1(t)}{\sqrt{f_2(t)}}
$$
and
\begin{equation}
\label{phidef}
\Phi(t) = \frac{1}{\pi} \cot^{-1} g(t).
\end{equation}
(Of course, $\Phi(t)=\A'(t)$ on $(r_1,r_2)$, but this new notation
will help avoid confusion soon.)  Then $g$ extends to a unique
holomorphic function on the (open) upper half plane.  The function
$\Phi$ extends as well to a unique holomorphic function on the upper
half plane, together with all of $\R$ except the points $-\rhL$, $0$,
$r_1$, $r_2$, $1$, and $1+\rhR$.  To see why, notice that $g(t)^2+1$
has only four roots, in particular, simple roots at each of $-\rhL$,
$0$, $1$, and $1+\rhR$.  There is always a holomorphic branch of the
arccotangent of a holomorphic function on a simply-connected domain,
as long as that function does not take on the values $\pm i$; this
fact gives us the analytic continuation of $\Phi$.  Of course, $\cot
\pi\Phi(t) = g(t)$ for all $t$ in the upper half plane.

For real $t$ (except $-\rhL$, $0$, $r_1$, $r_2$, $1$, and $1+\rhR$),
define
$$
\Philim(t) = \lim_{s \rightarrow 0+} \Phi(t+is)=\Phi(t).
$$
(We will use this notation to distinguish between the function
$\Philim$ on the real line and the function $\Phi$ on the upper half
plane.)
In order to apply Theorem~\ref{titchthm}, we must
determine the
real and imaginary parts of $\Philim(t)$.  Outside of $(r_1,r_2)$,
$\Real \Philim(t)$ is piecewise constant (in particular, constant
between the points where $\Philim$ is undefined) since $\Philim'(t)$
is imaginary there.  The integrability of $\Philim'(t)$ at $t=r_1$ and
$t=r_2$ implies that $\Philim(t)$ is continuous there, which implies
that $\Real \Philim(t) = \A'(t)$ for all $t \in (0,1)$.

To determine the behavior of $\Real \Philim(t)$ for $t \not\in (0,1)$,
we just have to see how much it changes by at $-\rhL$, $0$, $1$, and
$\rhR$, since it is constant on $(-\infty,-\rhL)$, $(-\rhL,0)$,
$(1,1+\rhR)$, and $(1+\rhR,\infty)$.  To do so, notice that if $\Phi'$
has a pole with residue $r$ at a point on the real axis, then
$\Philim$ changes by $-r\pi i$ as one moves from the left of that point
to its right.  (To see this, integrate over a small semi-circle in the
upper half plane, centered at the point.)  If $g(u)=\pm i$, then
\begin{eqnarray*}
\Res_{t=u} \Phi'(t) &=& \frac1\pi \lim_{t \rightarrow u}
-\frac{(t-u)g'(t)}{1+g(t)^2}\\
&=&
\frac1\pi \lim_{t \rightarrow u}
-\frac{(t-u)g'(t)}{(g(t)+g(u))(g(t)-g(u))} \\
&=&
-\frac{1}{2\pi g(u)}.
\end{eqnarray*}
Therefore, if $g(u)=\pm i$, then $\Real \Philim$ changes by
$\pm\frac12$ from the left of $u$ to its right.

To determine the precise sign of $g(u)$ when $g(u)=\pm i$, we will
need to know how $1/\sqrt{f_2(t)}$ behaves when analytically continued
through the upper half plane.  We know that it is positive on
$(r_1,r_2)$.  If one analytically continues it along any path through
the upper half plane that starts in $(r_1,r_2)$ and ends on the real
axis to the left of $r_1$, then the result is a negative imaginary
number (i.e., one with argument $-\pi i/2$).  Similarly, if the path
ends to the right of $r_2$, then the result is a positive imaginary
number.

Thus, $g(-\rhL) = -i\, \sgn f_1(-\rhL)$.  In fact, $g(-\rhL)=i$,
because
$$
f_1(-\rhL) = -\lambda^2-\rhoLR\lambda-\rhL\rhR - \rhL-\rhL^2 < 0.
$$
It follows that $\Real \Philim$ increases by
$\frac12$ at $-\rhL$.  Similarly, $\Real \Philim$ decreases by
$\frac12$ at $1+\rhR$.

The analysis at $0$ and $1$ is slightly more subtle.  We have $g(0) =
-i\, \sgn f_1(0)$, and $g(1) = i\, \sgn f_1(1)$.  It turns out that
$\sgn f_1(0) = \sgn f_1(r_1)$ and $\sgn f_1(1) = \sgn f_1(r_2)$.
To prove this claim, we will deal with $\Imag \Philim$.  (The results
about $\Imag \Philim$ will be needed later, so this approach is
worthwhile even though one might wish for a direct proof.)

It is impossible for $\Imag \Philim(t)$ to vanish for $t \in (0,r_1)
\cup (r_2,1)$, since $\Real \Philim(t)$ is $0$ or $1$ for such $t$,
but $\Philim(t)$ cannot be $0$ or $1$ (since otherwise $g$ would have
a singularity at $t$, as one can see from \eqref{phidef}).  Thus, the
sign of $\Imag \Philim(t)$ is constant for $t$ in each of $(0,r_1)$
and $(r_2,1)$.

The imaginary part of the arccotangent does not depend on the branch
used (since the values of the arccotangent always differ by a multiple
of $\pi$).  To determine the sign of the imaginary part of
$\Philim(t)$, we will use the fact that for real $u$ with $|u| > 1$,
\begin{equation}
\label{arccotsign}
\sgn \Imag \cot^{-1}(ui) = -\sgn u.
\end{equation}
In order to apply this formula to $\Philim$, we need to determine the
sign of $g$ on the axis.  We determined above how $1/\sqrt{f_2(t)}$
behaves.  Since
$$
\Phi(t) = \frac{1}{\pi}\cot^{-1}
\left(\frac{f_1(t)}{\sqrt{f_2(t)}}\right),
$$
we find, by combining the facts about $1/\sqrt{f_2(t)}$ with
\eqref{arccotsign}, that $\sgn \Imag \Philim(t) = \sgn f_1(t)$ for $t
\in (0,r_1)$, and $\sgn \Imag \Philim(t) = -\sgn f_1(t)$ for $t \in
(r_2,1)$.

Because $\sgn \Imag \Philim$ is constant, we can deduce two important
facts.  First, we see that $\sgn f_1$ must be constant on $(0,r_1)$
and $(r_2,1)$.  Notice that in fact it is constant on $[0,r_1]$ and
$[r_2,1]$, because we showed earlier that it cannot vanish at one of
the endpoints of one of these intervals unless the interval consists
only of one point.  Second, we see that $\sgn \Imag \Philim(t) = \sgn
f_1(r_1)$ for $t \in (0,r_1)$ and $\sgn \Imag \Philim(t) = -\sgn
f_1(r_2)$ for $t \in (r_2,1)$.

Thus, having taken a short detour, we can see that $g(0) = -i\, \sgn
f_1(r_1)$, and $g(1) = i\, \sgn f_1(r_2)$.  It follows that $\Real
\Philim$ increases by $\frac12$ at $0$ iff $f_1(r_1) < 0$, and
decreases by $\frac12$ at $0$ iff $f_1(r_1) > 0$.  Notice that these
are exactly the conditions under which $\A'(0)$ is $1$ or $0$,
respectively.  Similarly, $\Real \Philim$ increases by $\frac12$ at
$1$ iff $f_1(r_2) > 0$, and decreases by $\frac12$ at $1$ iff
$f_1(r_2) < 0$, and these are exactly the conditions under which
$\A'(1)$ is $0$ or $1$, respectively.

The information that we have obtained determines $\Real \Philim$, and
in fact shows that
$$
\Real \Philim(t) =
\begin{cases}
0 & \textup{if } t < -\rhL,\\
\frac12 & \textup{if } -\rhL < t < 0,\\
\A'(t) & \textup{if } 0 < t < 1 \ (\textup{and } t \ne r_1, r_2),\\
\frac12 & \textup{if }1 < t < 1+\rhR, \textup{ and}\\
0 & \textup{if } 1+\rhR < t.\\
\end{cases}
$$
In other words, $\Real \Philim = \A' + J'$.

We would like to apply Theorem~\ref{titchthm} to conclude that the
Hilbert transform of $\A' + J'$ must be $\Imag \Philim$;  for the
hypotheses of the theorem to be satisfied, we
need some integral estimates, in particular that the integrals
$$
\intr |\Phi(r+is)|^2 \,dr
$$
exist for all $s > 0$ and are bounded.  To prove existence of the
integral, we use the estimate $\Phi(t)=O(t^{-1})$ as $t \rightarrow
\infty$.  To verify this estimate, notice that as $t \rightarrow
\infty$ in the upper half plane, $g(t) =
-(\rhoLR/2+\lambda)^{-1}it+O(1)$, and for such $t$ we have
$$
\Phi(t) = \frac{1}{\pi}\cot^{-1}(-(\rhoLR/2+\lambda)^{-1}it+O(1)) =
\frac{i(\frac{\rhoLR}{2}+\lambda)}{\pi t} + O(t^{-2}) + k,
$$
for some integer $k$ depending on $t$ and the branch of the
arccotangent.  For large $t$, continuity implies that $k$ must be
constant, and our knowledge of the behavior of $\Phi$ on the real axis
tells us that $k=0$.  It follows that $\Phi(t)=O(t^{-1})$, so the
integrals must converge.  To prove boundedness, we need only show that
the integrals remain bounded as $s \rightarrow 0$.  To see that they
do, notice that the limiting integrand has singularities, but they are
only logarithmic singularities (since $\Phi'$ has poles of
order~$1$ there), so it is still integrable.

We can now verify that $\A(1)=\lambda$.  (This fact, which
is necessary for $\A$ to be admissible, was stated earlier, but the
proof was postponed.)  To determine
$\A(1)$, we need to integrate $\A'$ from $0$ to $1$.  If ${\mathcal
C}$ denotes a semi-circle of radius $R$ centered at $0$, lying in the
upper half plane, and oriented clockwise, then for $R >
\max(\rhL,1+\rhR)$, Cauchy's theorem implies that
\begin{eqnarray*}
\Real \int_{\mathcal C} \Phi(t) \,dt & = &
\Real \int_{-R}^{R} \Philim(t) \,dt \\
& = & \int_{-\rhL}^{1+\rhR} \A'(t)+J'(t)\:dt \\
& = & \A(1) + \frac{\rhoLR}2.
\end{eqnarray*}
Since $\Phi(t) = \frac{i(\rhoLR/2+\lambda)}{\pi t}+ O(t^{-2})$,
\begin{eqnarray*}
\int_{\mathcal C} \Phi(t) \,dt & = & \int_{\mathcal C}
\frac{i(\frac{\rhoLR}{2}+\lambda)}{\pi t}\,dt + O(R^{-1}) \\
& = & (-\pi i)\frac{i(\frac{\rhoLR}{2}+\lambda)}{\pi}+O(R^{-1})\\
& = & \lambda+\frac{\rhoLR}{2}+O(R^{-1}).
\end{eqnarray*}
Hence, $\A(1)=\lambda$.

Now Theorem~\ref{titchthm}
tells us that because
$$
\Philim(t) = \lim_{s \rightarrow 0+} \Phi(t+is)
$$
(except where undefined),
the imaginary part of $\Philim$ is the Hilbert transform
of the real part, i.e.,
$$
\widetilde{\A'}+\widetilde{J'} = \Imag \Philim.
$$
To complete the proof of \eqref{pos}, we need more information about
how $\Imag \Philim(t)$ behaves for $t \in [0,1]$.  We know that $\Imag
\Philim(t) = 0$ for $t \in (r_1,r_2)$, by the definition of $\Phi$,
so for such $t$,
$$
(A(t)-\A(t))(\widetilde{\A'}(t)+\widetilde{J'}(t)) = 0.
$$
If we can ensure that
\begin{equation}
\label{goalhere}
(A(t)-\A(t))(\widetilde{\A'}(t)+\widetilde{J'}(t)) \ge 0
\end{equation}
for all $t \in (0,r_1) \cup (r_2, 1)$, then we will have proved
\eqref{pos}.

We will deal first with the sign of $A(t)-\A(t)$.  Recall that $\A'$
is constant on $(0,r_1)$, and is either $0$ or $1$ (assuming $r_1 >
0$).  Because of the Lipschitz condition $0 \le A' \le 1$, it follows
that either $A'(t)-\A'(t) \ge 0$ for all $t \in (0,r_1)$, or
$A'(t)-\A'(t) \le 0$ for all such $t$, according as $\A'$ is $0$ or
$1$ on that interval.  Integrating and using $A(0)=\A(0)=0$ implies
that $A(t)-\A(t) \ge 0$ for $t \in (0,r_1)$ in the first case (where
$\A'(0)=0$), and $A(t)-\A(t) \le 0$ in the second (where $\A'(0)=1$).
Similarly, if $\A'(1)=0$ then $A(t)-\A(t) \le 0$ for $t \in (r_2,1)$,
and if $\A'(1)=1$ then $A(t)-\A(t) \ge 0$ for $t \in (r_2,1)$.
Therefore, to prove \eqref{goalhere}, we need only prove the same
inequalities as here, with $A(t)-\A(t)$ replaced by $\Imag
\Philim(t)$.

We have already shown that $\sgn \Imag \Philim(t) = \sgn f_1(r_1)$ for
$t \in (0,r_1)$, and that $\sgn \Imag \Philim(t) = -\sgn f_1(r_2)$ for
$t \in (r_2,1)$.  We know as well that $\A'(0)=0$ if $f_1(r_1) > 0$
and $\A'(0)=1$ if $f_1(r_1) < 0$, and that similarly, $\A'(1) = 0$ if
$f_1(r_2) > 0$ and $\A'(1)=1$ if $f_1(r_2) < 0$.  (Note that the only
possible conditions under which $r_1$ or $r_2$ are roots of $f_1$ are
$r_1=0$ and $r_2=1$, respectively.)  These conditions, when combined
with those derived in the previous paragraph, give us what we need.
We conclude that for $t \in (0,r_1) \cup (r_2,1)$, we have $A(t)-\A(t)
\ge 0$ iff $\Imag \Philim(t) \ge 0$.  Since $\Imag \Philim(t) =
\widetilde{\A'}(t) + \widetilde{J'}(t)$, and is $0$ for $t \in
(r_1,r_2)$, we see that for all $t \in (0,1)$ (and trivially for $t
\not\in (0,1)$ since then $A(t)=\A(t)$),
$$
(A(t)-\A(t))(\widetilde{\A'}(t) + \widetilde{J'}(t)) \ge 0.
$$
This inequality implies \eqref{pos}, which completes our proof.

Thus, $\A$ is indeed the unique admissible function such that $\V(\A)$
is maximal.  We leave to the reader the task of checking that applying
this result to the specific functional arrived at in
\propref{functionalprop} leads to the explicit formula for
$\Pl_{\alpha,\beta,\gamma}$ given in \thmref{mainthm}.

\section{Conjectures and Open Questions}
\label{sec-conj}

The theorems we have proved do not answer all the natural questions
about the typical plane partition in a box, or about random lozenge
tilings of hexagons.

Given a location $(x,y)$ in the normalized $\alpha,\beta,\gamma$ hexagon, we
can ask whether the probability of finding a vertical lozenge near $(x,y)$ is
given by $\Pl_{\alpha,\beta,\gamma}(x,y)$.  \thmref{mainthm} tells us that
this is true if we average over all $(x,y)$ in some macroscopic region.
However, it is conceivable that there might be small-scale fluctuations in the
probabilities that would even out on a large scale.  We believe that that is
not the case.

\begin{conj}
\label{cep}
Let $V$ be any open set in the $\alpha, \beta, \gamma$ hexagon
containing the four points at which $\Pl_{\alpha,\beta,\gamma}$ is
discontinuous.  As $\scalefactor \rightarrow \infty$, the probability
of finding a vertical lozenge at normalized location $(x,y) \not\in V$
is given by $\Pl_{\alpha,\beta,\gamma}(x,y)+o(1)$, where the $o(1)$
error bound is uniform in $(x,y)$ for $(x,y) \not\in V$.
\end{conj}

There is numerical evidence that \conjref{cep} is true.  Also, the
analogous result for random domino tilings of Aztec diamonds has been
proved in \cite{CEP}, and it is not hard to prove that the local
statistics for the one-dimensional case described in
Section~\ref{sec-intro} do in fact converge to i.i.d.\ statistics, so
it is plausible that \conjref{cep} is true.  A similar result should
also hold for higher-order statistics (probabilities of finding
configurations of several lozenges); one can deduce explicit
hypothetical formulas for these probabilities from the theorems and
conjectures in \cite{CKP}.

Of course, we do not need to restrict our attention to hexagons, but
can look at tilings with lozenges of any regions.  Hexagons do seem to
involve the most elegant combinatorics and analysis, but one can prove
in general that almost all tilings of large regions have normalized
height functions that cluster around the solution to the problem of
maximizing a certain entropy functional; see
\cite{CKP} for the details.

We also conjecture an analogue of the arctic circle theorem of
Jockusch, Propp, and Shor.  (See \cite{JPS} for the original proof, or
\cite{CEP} for the proof of a stronger version on which our conjecture
is based.)
Define the {\it arctic region} of a lozenge tiling to be the set of
lozenges connected to the boundary by sequences of adjacent lozenges
of the same orientation (where a lozenge is said to be adjacent to
another lozenge, or to the boundary, if they share an edge).

\begin{conj}[Arctic Ellipse Conjecture]
\label{aec}
Fix $\varepsilon > 0$.  The probability that the boundary of the
arctic region is more than a distance $\varepsilon$ (in normalized
coordinates) from the inscribed ellipse is exponentially small in
the scaling factor $\scalefactor$.
\end{conj}

There are also several questions for which we do not even have
conjectural answers.

\begin{quest}
\label{guess}
Is there a simple way to derive the results of
Section~\ref{sec-typical} without having to guess any formulas?
\end{quest}

Such a method would be much more pleasant than our approach.  A good
test case would be the following open question.

\begin{quest}
\label{q-analogue}
Is there a $q$-analogue of \thmref{mainthm}?
\end{quest}

Of course, the non-trivial situation is when $q \rightarrow 1$ as
$\scalefactor \rightarrow \infty$ (although we do not know the precise
relationship between $q$ and $\scalefactor$ that will lead to
interesting limiting behavior).  There is a simple $q$-analogue of
MacMahon's enumeration of boxed plane partitions (also due to
MacMahon), and a $q$-analogue of \propref{prodformula} (which can be
established by a proof similar to that of \propref{prodformula}).  We
believe that the answer to \questref{q-analogue} is yes, and that the
same approach should work, but the guess work required is likely to be
tedious.  We hope that further development of these techniques will
someday let one answer such questions more easily.

\begin{quest}
\label{space}
Is there an analogue of \thmref{mainthm} for ``space partitions'' (the
natural generalization of plane partitions from the plane to space)?
\end{quest}

The answer to \questref{space} may well be yes, but it is extremely
unlikely that similar techniques can be used to prove it.  (For
example, no analogue of MacMahon's formula is known, and there is
no reason to believe that one exists.)

\section*{Appendix.  Converting the Sum to an Integral}

In Section~\ref{sec-setup}, we had to convert a sum to an integral.
The sum was
\begin{equation}
\label{oursum}
\sum_{1\le i < j \le \nleft+\ell+\nright} \log \left(
\frac{B^{-1}(j/\nbig)-B^{-1}(i/\nbig)}{j/\nbig-i/\nbig}
\right)\left(\frac{1}{\nbig}\right)^2,
\end{equation}
and we interpreted it as a Riemann sum for the double integral
\begin{equation}
\label{ourint}
\iint_{0 \le s < t \le \frac{\lambda+\rhoLR}{1+\rhoLR}} \log
\frac{C^{-1}(t)-C^{-1}(s)}{t-s} \,ds \,dt,
\end{equation}
for some function $C$ which we have not yet specified.  In
Section~\ref{sec-setup}, we claimed that the difference between the
sum and the integral is $o(1)$, and that $C$ can be chosen so that
$B'$ and $C'$ nowhere differ by more than $O(1/\nbig)$.  In this
appendix, we will define $C$ and justify these claims.

The main obstacle is that $(B^{-1})'$ can be quite large (infinite, in
fact, when $B'=0$).  We will now define a modification $C$ of $B$
designed to keep $(C^{-1})'$ from being too large.  We put $C(0) = 0$
and $C(b_i/\nbig)=i/\nbig$ for $1 \le i \le \nleft+\ell+\nright$ (so
that $C(1) = (\lambda+\rhoLR)/(1+\rhoLR)$).  Between these points, we
will define $C$ so that it is a continuous, strictly increasing,
piecewise linear function on $[0,1]$ such that $C'$ is constant on
intervals $\left({i}/{\nbig},(i+1)/{\nbig}\right)$, $C'$ is never
smaller than $1/\nbig^2$ or greater than~$1$, and $|C'-B'| =
O(1/\nbig)$.  There is no canonical way to do this; one way
that works is as follows.  If $b_{i+1} > b_i+1$, then set
$C((b_{i+1}-1)/\nbig) = i/\nbig + 1/\nbig^2$ (and similarly if $b_1 >
1$ set $C((b_1-1)/\nbig) = 1/\nbig^2$).  Then interpolate linearly to
define $C$ in between the points at which we have defined it so far.
Notice that changing $B$ to $C$ does not change the sum \eqref{oursum},
since $C^{-1}(i/\nbig) = b_i/\nbig = B^{-1}(i/\nbig)$.

To begin, for $0 \le s < t \le \frac{\lambda+\rhoLR}{1+\rhoLR}$,
we define
$$
f(s,t) = \log \frac{C^{-1}(t)-C^{-1}(s)}{t-s}.
$$
Because $(C^{-1})' = O(\nbig^2)$, the mean value theorem implies that
$f(s,t) = O(\log \nbig)$.

Consider small squares of side length $\nbig^{-1}$, with
their sides aligned with the $s$-~and $t$-axes and their upper right
corners at $(s,t) = (i/\nbig,j/\nbig)$, for $1 \le i < j \le
\nleft+\ell+\nright$.
Each square has area $\nbig^{-2}$, and $f$ is bounded by $O(\log \nbig)$,
so we can safely remove up to $o(\nbig^2/\log \nbig)$ squares {}from the
domain of integration and the corresponding terms {}from the sum
without changing either by more than $o(1)$.  We will do so in order
to restrict our attention to squares on which we can bound the partial
derivatives of $f$.

We first remove the $O(\nbig^{3/2})$ squares containing some point $(s,t)$
with $|s-t| \le \nbig^{-1/2}$.

Next, we remove all squares containing some point $(s,t)$
satisfying $(C^{-1})'(s) \ge \nbig^{1/3}$ or $(C^{-1})'(t) \ge
\nbig^{1/3}$.  We can check as follows that there are at most
$O(\nbig^{5/3})$ such squares.  Since $C^{-1}$ is increasing and
has range contained in $[0,1]$,
the set of all
$t$ with $(C^{-1})'(t) \ge
\nbig^{1/3}$ has measure $O(\nbig^{-1/3})$.  Also, as $(s,t)$ varies
over each small square, $(C^{-1})'(s)$ and $(C^{-1})'(t)$ are
constant.  Hence, in this step we are removing at most
$O(\nbig^{-1/3}/\nbig^{-2})=O(\nbig^{5/3})$ squares.

Thus, we can restrict our attention to squares containing only points
$(s,t)$ with $|s-t| \ge \nbig^{-1/2}$, $(C^{-1})'(s) \le \nbig^{1/3}$,
and $(C^{-1})'(t) \le \nbig^{1/3}$.

Now we will estimate the difference between the sum and the
integral.  If we can show that on each square, $f$ varies by at most
$o(1)$ (uniformly for all squares), then we will be done.
To determine how much $f$ can vary over a square of side length
$\nbig^{-1}$,
we compute
$$
\frac{\partial f}{\partial t} =
\frac{(C^{-1})'(t)}{C^{-1}(t)-C^{-1}(s)}-\frac{1}{t-s}.
$$
The second term has absolute value at most $\nbig^{1/2}$ (since by
assumption $|s-t| \ge \nbig^{-1/2}$).  To bound the first term, we start
with the denominator.  Because $C' \le 1$ everywhere, we have
$(C^{-1})' \ge 1$, so $|C^{-1}(t)-C^{-1}(s)| \ge |s-t|$, and thus
$$
\left|\frac{(C^{-1})'(t)}{C^{-1}(t)-C^{-1}(s)}\right|
\le \frac{|(C^{-1})'(t)|}{|s-t|}
\le (C^{-1})'(t)\nbig^{1/2}.
$$
Finally, since $(C^{-1})'(t) \le \nbig^{1/3}$, we have
$$
\frac{\partial f}{\partial t} = O(\nbig^{5/6}) = o(\nbig).
$$
The same holds for ${\partial f}/{\partial s}$, of course.

Therefore, over one of the small squares of side length $\nbig^{-1}$,
$f$ can vary by at most $o(1)$.  Thus,
the sum \eqref{oursum} differs from the
integral \eqref{ourint} by at most $o(1)$.

\section*{Acknowledgements}

We thank Tom Brennan, Persi Diaconis, and Noam Elkies for useful
discussions, the anonymous referee for pointing out an oversight in
our proof, and Valerie Samn for helping with the figures.


\begin{thebibliography}{EKLP}

\bibitem[BH]{BH} H.\ W.\ J.\ Bl{\"{o}}te and H.\ J.\ Hilhorst,
{\it Roughening transitions and the zero-temperature triangular
Ising antiferromagnet\/}, J.\ Phys.\ A {\bf 15} (1982),
L631--L637.

\bibitem[CS]{CS} L.\ Carlitz and R.\ Stanley, {\it Branchings and
partitions\/}, Proc.\ Amer.\ Math.\ Soc.\ {\bf
53} (1975), 246--249.

\bibitem[CEP]{CEP} H.\ Cohn, N.\ Elkies, and J.\ Propp, {\it Local
statistics for random domino tilings of the Aztec diamond\/},
Duke Math.\ J.\ {\bf 85} (1996), 117--166, arXiv:math.CO/0008243.

\bibitem[CKP]{CKP} H.\ Cohn, R.\ Kenyon, and J.\ Propp, {\it A
variational principle for domino tilings\/}, J.\ Amer.\ Math.\ Soc.\ {\bf 14}
(2001), 297-346, arXiv:math.CO/0008220.

\bibitem[DT]{DT} G.\ David and C.\ Tomei, {\it The problem of the
calissons\/}, Amer.\ Math.\ Monthly {\bf 96} (1989), 429--431.

\bibitem[EKLP]{EKLP} N.\ Elkies, G.\ Kuperberg, M.\ Larsen, and J.\
Propp, {\it Alternating sign matrices and domino tilings\/}, J.\ Algebraic
Combin.\ {\bf 1} (1992), 111--132 and 219--234.

\bibitem[GT]{GT} I.~M.\ Gelfand and M.~L.\ Tsetlin, {\it
Finite-dimensional representations of the group of unimodular
matrices\/} (Russian), Dokl.\ Akad.\ Nauk SSSR {\bf 71} (1950),
825--828 (reprinted in English in {\it Izrail M.\ Gelfand: Collected
Papers\/}, Springer-Verlag, Berlin, 1988, vol.\ 2, pp.\ 653--656).

\bibitem[JP]{JP} W.\ Jockusch and J.\ Propp, {\it Antisymmetric
monotone triangles and domino tilings of quartered Aztec diamonds\/},
to appear in the Journal of Algebraic Combinatorics.

\bibitem[JPS]{JPS} W.\ Jockusch, J.\ Propp, and P.\ Shor, {\it Random
domino tilings and the arctic circle theorem\/}, preprint, 1995,
arXiv:math.CO/9801068.

\bibitem[LS]{LS} B.~F.\ Logan and L.~A.\ Shepp, {\it A
variational problem for random Young tableaux\/}, Adv.\ in Math.\ {\bf 26}
(1977), 206--222.

\bibitem[M]{M} P.~A.\ MacMahon, {\it Combinatory Analysis\/},
Cambridge University Press, 1915--16 (reprinted by Chelsea Publishing
Company, New York, 1960).

\bibitem[PW]{PW} J.\ Propp and D.\ Wilson, {\it Exact sampling with
coupled Markov chains and applications to statistical mechanics\/},
Random Structures and Algorithms {\bf 9} (1996), 223--252.

\bibitem[R]{R} W.\ Rudin, {\it Real and Complex Analysis\/},
McGraw-Hill, New York, 1987.

\bibitem[Th]{Th} W.~P.\ Thurston, {\it Conway's tiling groups\/},
Amer.\ Math.\ Monthly {\bf 97} (1990), 757--773.

\bibitem[Ti]{T} E.~C.\ Titchmarsh, {\it Introduction to the Theory of
Fourier Integrals\/}, Oxford University Press, London, 1937.

\bibitem[VK1]{VK1} A.~M.\ Vershik and S.~V.\ Kerov, {\it Asymptotics
of the Plancherel measure of the symmetric group and the limiting form
of Young tables\/}, Soviet Math. Dokl. {\bf 18} (1977), 527--531.

\bibitem[VK2]{VK2} A.~M.\ Vershik and S.~V.\ Kerov, {\it Asymptotics of
the largest and the typical dimensions of irreducible representations
of a symmetric group\/}, Functional Anal.\ Appl.\ {\bf
19} (1985), 21--31.

\end{thebibliography}
\end{document}